\documentclass[a4paper,draft]{article}
\input{epsf.sty}

\newcommand\ds{\displaystyle}

\newtheorem{thm}{\noindent Theorem}[section]

\newtheorem{cor}{\noindent Corollary}[section]
\newtheorem{prop}{\noindent Proposition}[section]
\newtheorem{rem}{\noindent Remark}[section]

\newenvironment{proof}{
\begin{trivlist}
\item[\hspace{\labelsep}{\it\noindent Proof. }]}
 {{\hfill\rule{2mm}{2mm}} \end{trivlist} \par}
\begin{document}
\title{ {\bf On the Asymptotic Behavior of First Passage Time \\
Densities for Stationary Gaussian Processes\\ 
and Varying Boundaries}\thanks{This work has been performed within a 
joint cooperation agreement between 
Japan Science and Technology Corporation (JST) and Universit\`a di Napoli 
Federico II, under partial support by MIUR and INdAM (GNCS).}}            
\author
{E. Di Nardo$^{(1)},$ A.G. Nobile$^{(2)},$ E. Pirozzi$^{(3)}$ and L.M. Ricciardi$^{(3)}$}

\date{}
\maketitle

\par\noindent
\begin{center}
{\it (1)} Dipartimento di Matematica, Universit\`a della
Basilicata, Campus Macchia Romana, Potenza, Italy,
Email: dinardo@unibas.it\\
{\it (2)} Dipartimento di Matematica e Informatica, Universit\`a di
Salerno, Via Allende,
Baronissi (SA), Italy, Email: nobile@unisa.it\\
{\it (3)} Dipartimento di Matematica e Applicazioni, Universit\`a di
Napoli Federico II,
Via Cintia, Napoli, Italy, Email: \{enrica.pirozzi, luigi.ricciardi\}@unina.it
\end{center}

\begin{abstract} 
Making use of a Rice-like series expansion, for a class of stationary Gaussian 
processes the asymptotic behavior of the first passage time  probability density 
function  through certain time-varying boundaries, including periodic boundaries, 
is determined. Sufficient conditions are then given such that the density 
asymptotically exhibits an exponential behavior when the boundary is either 
asymptotically constant or asymptotically periodic.  
\end{abstract}
{\bf Keywords:}\hspace{0.2cm}{Exponential trends; Simulation; Damped oscillatory covariance}
\par\noindent
{\bf AMS 2000 subject classification}\hspace{0.2cm}{Primary: 60G15$\;$ Secondary:
60G10; 60G40}

\section{Introduction}
\setcounter{equation}{0}
First-passage-time (FPT) probability density functions (pdf's) through 
generally time-de\-pen\-dent boundaries play an essential role in many 
applied fields including the stochastic description of the behavior 
of various biological systems (see, for instance, \cite{Dicrescenzo00}, 
\cite{Kosty81},  \cite{Lansky99}, \cite{Ricciardi77},  \cite{Ricciardi99} and the references
therein). Investigations have  essentially proceeded along the following three main
directions:  {\em (i)\/} to search for  closed-form solutions under suitable assumptions 
on the considered  stochastic processes and on the boundaries (see, for instance, 
\cite{Dicrescenzo97}, \cite{Dinardo01}, \cite{Giorno88}, \cite{Giorno89two}, \cite{Gutierrez97}); 
{\em (ii)\/} to devise  numerical algorithms to evaluate FPT densities (see, for 
instance, \cite{Anderssen73},\cite{Buonocore87}, \cite{Buonocore90}, \cite{Daniels96}, \cite{Daniels00}, \cite{Durb71}, \cite{Durb85}, \cite{Ferebee82}, \cite{Giorno89}, \cite{Roberts68})  and {\em(iii)\/} to  analyze the 
asymptotic behavior of the FPT
densities as  boundaries or time grow larger (see, for instance, \cite{Giorno90}, 
\cite{Nobile85}, \cite{Nobile85two}, \cite{Sacerdote88}, \cite{Sato77}). 
The present paper, that falls within category {\em (iii)}, 
is the natural extension of previous investigations 
carried out by us for the class of one-dimensional diffusion processes 
admitting steady state densities in the presence of single asymptotically 
constant boundaries or of single asymptotically periodic boundaries
(\cite{Giorno90}, \cite{Nobile85}, \cite{Nobile85two}). In such cases, computational as well as
analytical  results have indicated that the FPT pdf through an asymptotically periodic 
boundary is susceptible of an excellent non-homogeneous exponential 
approximation for large times and for large boundaries.
\par
However, if one deals with problems involving processes characterized by 
memory effects, or evolving on a time scale 
which is comparable with that of measurements or observations, 
the customarily assumed strong Markov property does not hold any longer; 
hence facing FPT problems for non Markovian 
processes becomes unavoidable. As is well known, for  
such processes no manageable equation holds for the conditional FPT pdf:   
only an excessively cumbersome series expansion  is available  
when the process is Gaussian,  stationary and mean square differentiable 
(cf. \cite{Ricciardi83}, \cite{Ricciardi86} and the references therein).  
\par
Due to the outrageous  complexity exhibited by the numerical evaluation of 
the involved partial sums on accounts of the analytical form of the involved 
terms,  a totally  different approach has been recently undertaken  in order 
to obtain information on the asymptotic behavior of the FPT densities for 
a class of normal processes. This consists of 
a simulation procedure  (see \cite{Dinardo97} and \cite{Dinardo00}) implemented  
to generate sample paths and to estimate the corresponding 
FPT densities.  Extensive 
computations have thus been performed  to gain some insight on the behavior of  
the FPT pdf  through varying boundaries. The results of the simulations, obtained by 
means of  a parallel supercomputer CRAY T3E,  
have indicated that for certain periodic boundaries not very distant from the initial 
value of the process, the simulated FPT pdf, $\tilde{g}(t)$, soon exhibits damped 
oscillations having the same period of the boundary. Indeed, to a  high degree of 
accuracy, $\tilde{g}(t)$  can be represented in the form
\begin{eqnarray}
\tilde{g}(t)\sim {\tilde \beta}(t)\,e^{-{\widehat\alpha}\,t},
\label{eq:(1.1)}
\end{eqnarray}
with $\widehat\alpha$ and $\tilde\beta(t)$ specified by means of the data 
obtained via the performed simulations \cite{Dinardo01two}. Note  that 
(\ref{eq:(1.1)}) can  be thrown in the equivalent form
\begin{eqnarray}
\tilde{g}(t)\sim {\tilde a}(t)\,\exp\Bigl\{-\int_0^t {\tilde a}(\tau)\,d\tau\Bigr\},
\label{eq:(1.2)}
\end{eqnarray}
where ${\tilde a}(t)>0$ is a periodic function having the same period of the boundary. 
Hence,  for periodic boundaries, even though 
not very distant from the initial position of the process, the estimated FPT pdf 
 appears to admit a non-homogeneous exponential approximation.  
\par
In the present paper the relevance and the validity of such an unexpected numerical 
result is confirmed. Indeed, it will be proved analytically that the 
non-homogeneous  exponential approximation (\ref{eq:(1.2)})  holds for a wide class 
of stationary Gaussian processes in the presence of boundaries that either 
possess a horizontal asymptote or are asymptotically periodic. 
\par
In Section 2 we shall briefly recall some basic notation 
that will be used throughout this paper; in Section 3 we shall assume 
that  the boundaries  possess a horizontal asymptote, and in Section 4 
that they are asymptotically periodic.  Finally, in Section 5 
for a stationary Gaussian process with zero mean and damped oscillatory covariance, 
the simulated FPT pdf $\tilde{g}(t)$ is compared with  
the non-homogeneous exponential approximation  for the FPT pdf. 
%
\section{Mathematical background}
\setcounter{equation}{0}
Let $\left\{ X(t),t\geq 0\right\} $ be a one-dimensional, non-singular 
stationary Gaussian process with mean $E[X(t)]=0$ and covariance  
$E[X(t)X(\tau )]=\gamma (t-\tau )$ such that 
$\gamma (0)=1,\dot{\gamma}(0)=0$ and $\ddot{\gamma}(0)<0$. Then $\dot{X}(t)$, 
the derivative of $X(t)$ with respect to $t$, exists in the mean-square sense. 
Let $S(t)\in C^1[0,+\infty )$ be an arbitrary function 
such that $X(0)=x_0<S(0)$. Then, 
\begin{eqnarray}
T=\inf_{t\geq 0}\{t:X(t)>S(t)\},\quad X(0)=x_0
\label{eq:FPTrandomvariable}
\end{eqnarray}
is the FPT random variable and
\begin{eqnarray}
g(t|x_0)={\partial\over\partial t} P(T<t)
\label{eq:FPTdensity}
\end{eqnarray}
is the FPT pdf of $X(t)$ through $S(t)$ conditional upon $X(0)=x_0$. 
For all $n\in {\bf N}$ and $0=t_0<t_1<\ldots <t_n$ we denote by 
$W_n(t_1,\ldots ,t_n|x_0)\;dt_1\cdots dt_n$ the probability that $X(t)$ 
crosses $S(t)$ from below in the intervals $(t_1,t_1+dt_1),\ldots, 
(t_n,t_n+dt_n)$ given that $X(0)=x_0$. 
As shown in \cite{Ricciardi83}, the functions $W_n$ can be expressed
as
\begin{eqnarray}
&& \hspace*{-2.0cm} W_n(t_1,\ldots ,t_n|x_0)=\int_{\dot{S}(t_1)}^{+\infty} dy_1
\int_{\dot{S}(t_2)}^{+\infty} dy_2\cdots 
\int_{\dot{S}(t_n)}^{+\infty}  \, \prod_{i=1}^n\,[y_i-\dot{S}(t_i)] \nonumber\\ 
&&\hspace*{1.2cm} \times p_{2n}[S(t_1),t_1;\ldots ;S(t_n),t_n;y_1,t_1;\ldots ;
y_n,t_n|x_0]\;dy_n,
\label{eq:def1W_n}
\end{eqnarray}
where $p_{2n}(x_1,t_1;\ldots ,x_n,t_n;y_1,t_1;\ldots ;y_n,t_n|x_0)$ is the joint 
pdf of $2n$ random variables $X(t_1),$ $\ldots,X(t_n),$ $Y(t_1)=\dot{X}(t_1),\ldots$, 
$Y(t_n)=\dot{X}(t_n)$ conditional upon $X(0)=x_0$:
\begin{eqnarray}
& &\hspace*{-0.7cm} p_{2n}(x_1,t_1;\ldots ;x_n,t_n; y_1,t_1;\ldots ;y_n,t_n|x_0)
= {1\over (2\pi)^n 
\left| \Lambda_{2n+1}(t_1,\ldots ,t_n) \right|^{1/2}}
\label{eq:defp_2n}\\ 
& &\hspace*{-0.3cm} \times \exp\biggl\{-{1\over 2 \left| \Lambda_{2n+1}
(t_1,\ldots ,t_n) \right|} \sum_{i,j=1}^{2n} 
l_{i+1,j+1}(t_1,\ldots, t_n)\, \bigl[\widehat{x}_i-x_0\,\widehat{\gamma}(t_i)\bigr]\,
\bigl[\widehat{x}_j-x_0\,\widehat{\gamma}(t_j)\bigr]\biggr\}.\nonumber
\end{eqnarray}
Here $l_{i+1,j+1}(t_1,\ldots ,t_n)$ denotes the cofactor of the element 
$\lambda_{i+1,j+1}(t_1,\ldots ,t_n)$ of the covariance matrix 
$\Lambda_{2n+1}(t_1,\ldots,t_n)$ of 
$X(0),X\left( t_1\right) ,...,X\left( t_n\right) ,
\dot{X}\left( t_1\right) ,...,\dot{X}\left( t_n\right)$, i.e.
\begin{eqnarray}
\lambda_{i+1,j+1}(t_1,\ldots ,t_n) = \left\{ \begin{array}{l}
E[X(t_i)X(t_j)] = \gamma(t_i-t_j) = \gamma(t_j-t_i), \\
\hspace*{3.12cm} (i=0,1,\ldots, n, j=0,1,\ldots, n) \\
\\
E[X(t_i)\dot{X}(t_{j-n})] = - \dot{\gamma}(t_i-t_{j-n}) =
\dot{\gamma}(t_{j-n}-t_i), \\
\hspace*{3.12cm} (i=0,1,\ldots, n, j=n+1, \ldots, 2n) \\
\\
E[\dot{X}(t_{i-n})X(t_j)] = \dot{\gamma}(t_{i-n}-t_j) = - 
\dot{\gamma}(t_j-t_{i-n}), \\
\hspace*{3.12cm} (i=n+1,\ldots, 2n, j=0, \ldots, n) \\
\\
E[\dot{X}(t_{i-n})\dot{X}(t_{j-n})] = - \ddot{\gamma}(t_{i-n}-t_{j-n}) 
= - \ddot{\gamma}(t_{j-n}-t_{i-n}), \\
\hspace*{3.12cm} (i=n+1\ldots, 2n, j=n+1\ldots, 2n),
\end{array} \right.
\label{eq:lambda}
\end{eqnarray}
\begin{eqnarray*}
&&\widehat{x}_i:=\left\{ \begin{array}{ll}
x_i,&i=1,2,\ldots,n\\
y_{i-n},&i=n+1,n+2,\ldots,2n
\end{array} \right.\\
&&\\
&& \widehat{\gamma}(t_i):=\left\{ \begin{array}{ll}
\gamma(t_i),&i=1,2,\ldots,n\\
\dot{\gamma}(t_{i-n}),&i=n+1,n+2,\ldots,2n
\end{array} \right.
\end{eqnarray*}
and $|A|$ denotes the determinant of a matrix $A$. 
Substituting (\ref{eq:defp_2n}) in (\ref{eq:def1W_n})  one has:
\begin{eqnarray}
& &\hspace*{-0.5cm} W_n(t_1,\ldots ,t_n|x_0)=  {1\over (2\pi)^n 
\left| \Lambda_{2n+1}(t_1,\ldots ,t_n) \right|^{1/2}}
\int_{\dot{\psi}(t_1)}^{+\infty} d\xi_1 \int_{\dot{\psi}(t_2)}^{+\infty} d\xi_2\cdots 
\label{eq:def2W_n}\\
& &\hspace*{1.0cm}\cdots\,\int_{\dot{\psi}(t_n)}^{+\infty}
\prod_{i=1}^n\,\bigl[\xi_i-\dot{\psi}(t_i)\bigr]\, 
H(t_1,\ldots ,t_n;\xi_1,\ldots ,\xi_n)\, 
K(t_1,\ldots ,t_n;\xi_1,\ldots ,\xi_n)\;d\xi_n\, ,\nonumber
\end{eqnarray}
where
\begin{eqnarray}
& &\hspace*{-0.7cm} \psi(t_i)  =  S\left(t_i\right) - x_0\,\gamma\left(t_i\right)
\qquad (i=1,2,\ldots,n)\nonumber\\ 
& &\hspace*{-0.7cm} H(t_1,\ldots ,t_n;\xi_1,\ldots ,\xi_n)=\exp\biggl\{ 
-{1\over 2 \left| 
\Lambda_{2n+1}(t_1,\ldots ,t_n) \right|} \sum_{i,j=1}^n l_{i+n+1,j+n+1}(t_1,
\ldots, t_n)\,\xi_i\,\xi_j\biggr\}\nonumber\\
& &\label{eq:componentW_n}\\
& &\hspace*{-0.7cm} K(t_1,\ldots ,t_n;\xi_1,\ldots ,\xi_n)
=\exp\biggl\{ -{1\over 2 \left| 
\Lambda_{2n+1}(t_1,\ldots ,t_n) \right|}\, \biggl[\,\sum_{i,j=1}^n l_{i+1,j+1}(t_1,
\ldots, t_n)\, \psi(t_i)\, \psi(t_j)\nonumber\\
& &\hspace*{1.5cm}+\sum_{i,j=1}^n l_{i+1,j+n+1}(t_1,\ldots, t_n)\, \psi(t_i)\,\xi_j
+\sum_{i,j=1}^n l_{i+n+1,j+1}(t_1,\ldots, t_n)\,\psi(t_j)\,\xi_i\biggr]
\biggr\}.\nonumber
\end{eqnarray}
As shown in \cite{Ricciardi83},  $g(t|x_0)$ can be expressed 
as the following  Rice-like series: 
\begin{eqnarray}
g(t|x_0)=W_1(t|x_0)+\sum_{i=1}^\infty
(-1)^i\int_0^tdt_1\int_{t_1}^tdt_2\cdots
\int_{t_{i-1}}^t W_{i+1}(t_1,\ldots ,t_i,t|x_0)\;dt_i\,, 
\label{eq:FPTGaussiandensity} 
\end{eqnarray}
with $x_0< S(0)$.
\par 
Let now $a_r$ denote the partial sum of order $r$ of the series in 
(\ref{eq:FPTGaussiandensity}). Then, for each $t>0$ the partial sums  of even order 
give a lower bound to $g$, whereas the partial sums  of odd order 
provide an upper bound to $g$. Since the evaluation of 
the partial sums is very cumbersome because of the complexity of the 
functions $W_n$ and of their integrals, a first
approximation of FPT density can be carried out by evaluating $W_1(t|x_0)$. 
The explicit expression of $W_1(t|x_0)$ (cf. \cite{Ricciardi86}) is:
\begin{eqnarray}
& & W_1(t|x_0)  =  {|\Lambda_3(t)|^{1/2}\over 2\,\pi\, [1-\gamma^2(t)]}
\exp \left\{ -{[S(t)-x_0\,\gamma (t)]^2\over 2\,[1-\gamma^2(t)]}\right\} 
\label{eq:W_1}\\
& &\hspace*{1.5cm} \times  \left[ \exp \left\{ -{\sigma^2(t|x_0)\over 2}\right\} 
-\sqrt{{\pi\over 2}}\,\sigma(t|x_0)\,\hbox{Erfc} \left( 
{\sigma(t|x_0)\over \sqrt{2}} \right) \right],\nonumber
\end{eqnarray}
where
\begin{eqnarray}
& &|\Lambda_3 (t)|=-\ddot{\gamma}(0)\,\bigl[1-\gamma^2(t)\bigr]
-\bigl[\dot{\gamma}(t)\bigr]^2\nonumber\\
& &\label{eq:elementW_1}\\
& &\sigma (t|x_0)=\left( {1-\gamma^2(t)\over |\Lambda_3 (t)|}\right)^{1/2}
\left\{ \dot{S}(t)+{\dot{\gamma}(t)\,\bigl[\gamma (t)\,S(t)-x_0\bigr]
\over 1-\gamma^2(t)}\right\}\nonumber
\end{eqnarray}
and
\begin{eqnarray}
\hbox{Erfc}(z) = {2\over\sqrt{\pi}}\,\int_{z}^{+\infty} 
exp\Bigl\{-{y^2\over 2}\Bigr\}\;dy, \qquad z \in {\bf R}.
\label{eq:erfc}
\end{eqnarray}
\par
We stress that although (\ref{eq:FPTGaussiandensity}) gives 
a formal analytical expression for the FPT densities 
through arbitrary time-dependent boundaries, no reliable numerical evaluations 
appear to be feasible due to the  complexity of (\ref{eq:def2W_n}) 
and (\ref{eq:componentW_n}). 
Furthermore, for all $t>0$ the first-order approximation $W_1(t|x_0)$, that 
provides an upper bound to the FPT pdf in (\ref{eq:FPTGaussiandensity}), 
is a good approximation of $g$ only for small values of $t$. 
%
\section{Asymptotically constant boundary}
\setcounter{equation}{0}
In this Section we consider the  FPT  problem for an asymptotically 
constant boundary
\begin{eqnarray}
S(t) = S_0 + \varrho(t), \quad t\geq 0,
\label{eq:asymconstboundary} 
\end{eqnarray}
with $S_0 \in {\bf R}$ and where $\varrho(t) \in C^1[0,+\infty)$ is a bounded function  
independent of $S_0$ and such that 
\begin{eqnarray}
\lim_{t \to +\infty} \varrho(t) = 0 \quad \hbox{and} \quad 
\lim_{t \to +\infty} \dot{\varrho}(t) = 0.
\label{eq:asymconstantrho}
\end{eqnarray}
The following proposition shows that under suitable hypotheses   
on the covariance $\gamma(t)$, the function $W_1(t|x_0)$, given in 
(\ref{eq:W_1}), approaches a constant value as $t$ increases.  
%
\begin{prop}{\em 
If 
\begin{eqnarray}
\lim_{t \to +\infty} \gamma(t) = 0 \quad \hbox{and} \quad 
\lim_{t \to +\infty} \dot{\gamma}(t) = 0,
\label{eq:asymconstantgamma1}
\end{eqnarray}
then
\begin{eqnarray}
R(S_0):=\lim_{t \to +\infty} W_1(t|x_0) =
{\sqrt{-\ddot{\gamma}(0)}\over 2\,\pi}\, \exp \left\{ - {S_0^2\over 2} \right\}.
\label{eq:R(S_0)}
\end{eqnarray}
}\end{prop}
\begin{proof}
Due to (\ref{eq:asymconstantrho}), from (\ref{eq:asymconstboundary}) it follows
\begin{eqnarray}
\lim_{t \to +\infty} S(t) = S_0, \qquad\qquad
\lim_{t \to +\infty} \dot{S}(t) = 0.
\label{eq:asymconstboundary1}
\end{eqnarray}
Furthermore, by virtue of (\ref{eq:asymconstantgamma1}) and 
(\ref{eq:asymconstboundary1}), from (\ref{eq:elementW_1}) one has
\begin{eqnarray}
\lim_{t \to +\infty}|\Lambda_3(t)| = - \ddot{\gamma}(0),\qquad\qquad
\lim_{t \to +\infty} \sigma(t|x_0) = 0.
\label{eq:asymconstcondition}
\end{eqnarray}
Taking the limit as $t\to +\infty$ in (\ref{eq:W_1}), and making use of 
(\ref{eq:asymconstboundary1}) 
and (\ref{eq:asymconstcondition}), we are finally led to (\ref{eq:R(S_0)}). 
This completes the proof. 
\end{proof}
We note that $R(S_0) > 0$ for all $S_0 \in {\bf R}$ and
\begin{eqnarray}
\lim_{S_0 \to +\infty} R(S_0) = 0. 
\label{eq:asymR(S_0)}
\end{eqnarray}
The nonzero asymptotic value given by (\ref{eq:R(S_0)}), {\em per se} indicates 
the inadequacy of $W_1(t|x_0)$ to provide a valid approximation to $g(t|x_0)$ for large 
time. On the contrary, the goodness of such an approximation for small times is 
confirmed. 
%
\begin{thm}\label{thmconst}{\em
Let $S(t)$, $t\geq 0$, be bounded and such that (\ref{eq:asymconstboundary}) 
and (\ref{eq:asymconstantrho}) hold. If 
\begin{eqnarray}
\lim_{t \to +\infty} \gamma(t) = 0,\qquad 
\lim_{t \to +\infty} \dot{\gamma}(t) = 0,\qquad 
\lim_{t \to +\infty} \ddot{\gamma}(t) = 0
\label{eq:asymconstantgamma2}
\end{eqnarray}
and
\begin{eqnarray}
\lim_{S_0 \to +\infty} {\varrho\left({\ds{t\over R(S_0)}}\right)\over S_0} = 0,
\label{eq:hypthmconst}
\end{eqnarray}
with $R(S_0)$ defined in (\ref{eq:R(S_0)}), then
\begin{eqnarray}
\lim_{S_0 \to +\infty} {1\over R(S_0)} \; g\left({t\over R(S_0)} \;  
\Big| \;  x_0\right) = e^{-t}, \qquad x_0 < S(0).
\label{eq:constant_density}
\end{eqnarray}
}\end{thm}
\begin{proof}
Since $R(S_0)>0$ for all $S_0 \in {\bf R}$, changing $t$ in $t/R(S_0)$ 
in (\ref{eq:FPTGaussiandensity}) we obtain:
\begin{eqnarray*}
&&\hspace*{-0.6cm} {1\over R(S_0)} \, g\left( {t\over R(S_0)}\;  \Big| \; x_0\right) 
= {1\over R(S_0)} W_1 \left({t\over R(S_0)}
 \;  \Big| \; x_0\right) \\ 
& & \hspace*{-0.5cm} + \, \, {1\over R(S_0)} \sum_{i=1}^\infty (-1)^i 
\int_0^{{t\over R(S_0)}} 
dt_1 \int_{t_1}^{{t\over R(S_0)}} dt_2 \cdots \int_{t_{i-1}}^{{t\over R(S_0)}} 
\,\, W_{i+1}\left(t_1,\ldots ,t_i,{{t\over R(S_0)}} \,\,  \Big| \,\,x_0 \right)\; dt_i 
\end{eqnarray*}
and hence,
\begin{eqnarray}
& & \hspace*{-0.6cm} {1\over R(S_0)} \, g\left( {t\over R(S_0)} \,\,  
\Big| \,\,x_0\right) = {1\over R(S_0)} 
W_1 \left({t\over R(S_0)} \;  \Big| \; x_0\right) 
\label{eq:thmconstant2} \\
&  & \hspace*{-0.5cm} + \sum_{i=1}^\infty  {(-1)^i\over [R(S_0)]^{i+1}} 
\int_0^{t} d\tau_1 \int_{\tau_1}^{t} d\tau_2 \cdots \int_{\tau_{i-1}}^{t} 
\, W_{i+1}\left( {\tau_1\over R(S_0)},\ldots ,{\tau_i\over R(S_0)},
{t\over R(S_0)} \,\,  \Big| \,\, x_0 \right)\; d\tau_i \nonumber
\end{eqnarray}
after the change of variables $\tau_i = t_i\,R(S_0)$ for $i=1,2, \ldots$. 
\par
We now prove that for all $n=1,2,\ldots$ one has:
\begin{eqnarray}
& &\lim_{S_0 \to +\infty} {1\over \bigl[R(S_0)\bigr]^n} \, W_n \left( 
{\tau_1\over R(S_0)},\ldots,{\tau_n\over R(S_0)} \;  \Big| \;x_0 \right) = 1, 
\label{eq:thmconstant3} \\
& & \hspace*{4cm}(0 < \tau_1 < \tau_2 < \ldots < \tau_n). \nonumber
\end{eqnarray}
To simplify the notation, we set: 
\begin{eqnarray}
\vartheta_1={\tau_1\over R(S_0)},\ldots, \vartheta_n={\tau_n\over R(S_0)}\qquad
(0 < \vartheta_1 < \vartheta_2 < \ldots < \vartheta_n).
\label{eq:thmconstant3bis}
\end{eqnarray}
Hence, recalling (\ref{eq:def2W_n}), for $n=1,2,\ldots$ one has:
\begin{eqnarray}
& & {1\over [R(S_0)]^n} W_n \left( \vartheta_1,\ldots,
\vartheta_n \;\Big| \; x_0 \right)=
{ 1\over (2\pi)^n \left| \Lambda_{2n+1} \left( \vartheta_1,\ldots,
\vartheta_n \right) \right|^{1/2} } 
\label{eq:thmconstant4}\\
& &\hspace*{1.0cm} \times \int_{\dot{\psi}\left(\vartheta_1\right)}^{+\infty} d\xi_1
\int_{\dot{\psi}\left(\vartheta_2\right)}^{+\infty} d\xi_2
\cdots \int_{\dot{\psi}\left(\vartheta_n\right)}^{+\infty} 
\prod_{i=1}^n \left[ \xi_i-\dot{\psi}\left( \vartheta_i\right) \right]\,
H \left( \vartheta_1,\ldots,\vartheta_n
;\xi_1,\ldots,\xi_n \right)\,\nonumber\\
& & \hspace*{1.0cm} \times  \Biggl\{{1\over \bigl[R(S_0)\bigr]^n}\,
K\left( \vartheta_1,\ldots,\vartheta_n;
\xi_1,\ldots,\xi_n \right)\Biggr\}\;d\xi_n\, ,\nonumber 
\end{eqnarray}
where the functions $\psi$, $H$  and $K$ are defined in  (\ref{eq:componentW_n}). 
\par
We shall now make use of the following relations (see Appendix I):
\begin{eqnarray}
& & \lim_{S_0 \to +\infty} \dot{\psi}\left( \vartheta_i\right) =0\qquad(i=1,2,\ldots,n),
\label{eq:thmconstant10}\\
& & \lim_{S_0 \to +\infty} \, \Big| \, \Lambda_{2n+1} \left(\vartheta_1, 
\ldots,\vartheta_n \right) \, \Big| \, = \bigl[-\ddot{\gamma}(0)\bigr]^n,
\label{eq:thmconstant7}\\
& & \lim_{S_0 \to +\infty} 
H\left( \vartheta_1,\ldots,\vartheta_n;\xi_1,\ldots,\xi_n \right)
=\exp \left\{ - {1\over 2\,[-\ddot{\gamma}(0)]} \sum_{i=1}^n
\xi_i^2 \right\},\label{eq:thmconstant9}\\
& & \lim_{S_0 \to +\infty} \Biggl\{{1\over \bigl[R(S_0)\bigr]^n}\,
K\left( \vartheta_1,\ldots,\vartheta_n;\xi_1,\ldots,\xi_n \right)\Biggr\} 
={(2\pi)^n\over[- \ddot{\gamma}(0)]^{n/2}}\,\cdot
\label{eq:thmconstant11ter}
\end{eqnarray}
Taking the limit as $S_0 \to +\infty$ in (\ref{eq:thmconstant4}) and 
recalling (\ref{eq:thmconstant10}), (\ref{eq:thmconstant7}),  
(\ref{eq:thmconstant9}) and  (\ref{eq:thmconstant11ter}) one 
then obtains: 
\begin{eqnarray}
& & \lim_{S_0 \to +\infty} {1\over [R(S_0)]^n} W_n \left(\vartheta_1,\ldots,
\vartheta_n \,\, \big| \,\, x_0 \right) 
\label{eq:thmconstant12} \\
& & \hspace*{1cm}  = {1\over [-\ddot{\gamma}(0)]^n}  \, \int_0^{+\infty} 
d \xi_1 \int_0^{+\infty} d \xi_2 \cdots \int_0^{+\infty} \exp \left\{ - 
{1\over 2\,[-\ddot{\gamma}(0)]} \sum_{i=1}^n \xi_i^2 \right\}
\prod_{i=1}^n \xi_i\, \;d \xi_n  
\nonumber \\
& & \hspace*{1cm} = {1\over [-\ddot{\gamma}(0)]^n} \left[ \int_0^{+\infty} \xi\, 
\exp \left\{ - {\xi^2\over 2\,[-\ddot{\gamma}(0)]}  \right\}\; d \xi \right]^n.\nonumber
\end{eqnarray}
Since 
$$
\int_0^{+\infty} \xi\, \exp \left\{ - {\xi^2\over 2\,[-\ddot{\gamma}(0)]} 
 \right\} d \xi  = 2\, [-\ddot{\gamma}(0)]\,\int_0^{+\infty} z\, e^{- z^2} dz 
=-\ddot{\gamma}(0),
$$
relation (\ref{eq:thmconstant3}) immediately follows from 
(\ref{eq:thmconstant12}). Due to 
(\ref{eq:thmconstant3}), taking the limit as $S_0 \to +\infty$ in 
(\ref{eq:thmconstant2}), one finally obtains:
\begin{eqnarray}
\lim_{S_0 \to +\infty} {1\over R(S_0)}\; g\left( {t\over R(S_0)} \; \Big| \; x_0\right) 
= 1 + \sum_{i=1}^\infty (-1)^i 
\int_0^{t} d\tau_1 \int_{\tau_1}^{t} d\tau_2 \cdots \int_{\tau_{i-1}}^{t} 
d\tau_i\,, 
\label{eq:thmconstant13}
\end{eqnarray}
that identifies with (\ref{eq:constant_density}). The proof is thus complete. 
\end{proof}
The following corollary is an immediate consequence of Theorem \ref{thmconst}.
\begin{cor}{\em
Under the assumptions of Theorem \ref{thmconst}, for $S_0\to +\infty$ one 
has:
\begin{eqnarray}
g(t \,\vline\,  x_0) \sim R(S_0)\,\exp\{-R(S_0)\,t\},
\qquad\qquad\forall t>0,
\label{eq:comp_g(t)}
\end{eqnarray}
with $R(S_0)$ defined in (\ref{eq:R(S_0)}). 
}\end{cor}
This Corollary expresses the asymptotic exponential trend of the FPT density 
as the boundary moves away from the process'  starting point.
\section{Asymptotically periodic boundary}
\setcounter{equation}{0}
The  FPT  problem  in the case of an asymptotically periodic boundary 
will be the object of the present Section. More specifically, we shall 
focus our attention on boundaries of the form
\begin{eqnarray}
S(t) = S_0 + \varrho(t), \quad  t\geq 0,
\label{eq:asymperiodicboundary} 
\end{eqnarray}
where $S_0 \in {\bf R}$ and $\varrho(t) \in C^1[0,+\infty)$ is a bounded function  
independent of $S_0$  and such that
\begin{eqnarray}
\lim_{k \to \infty} \varrho(t+k\,Q) = Z(t),\qquad 
\lim_{k \to \infty} \dot{\varrho}(t+k\,Q) = \dot{Z}(t),
\label{eq:Z(t)}
\end{eqnarray}
where $Z(t)$  is a periodic function of period $Q > 0$  satisfying
\begin{eqnarray}
\int_0^Q Z(\tau)\;d\tau = 0.
\label{eq:intZ(t)}
\end{eqnarray}
%
\begin{prop} {\em
If (\ref{eq:asymconstantgamma1}) holds, then 
\begin{eqnarray}
& &R[Z(t)]:=\lim_{k \to \infty} W_1(t+k\,Q\,|\,x_0) =
{\sqrt{-\ddot{\gamma}(0)}\over 2\,\pi} \exp \left\{ - {[S_0+Z(t)]^2\over 2} \right\} 
\label{eq:(4.27)}\\
& &\hspace*{2cm} \times  \left[ \exp \left( - {[\dot{Z}(t)]^2
\over 2\,[-\ddot{\gamma}(0)]} 
\right)  -  \sqrt{{\pi\over 2\,[-\ddot{\gamma}(0)]}} \, \dot{Z}(t) \, {\rm Erfc} 
\, \left({\dot{Z}(t)\over \sqrt{2\,[-\ddot{\gamma}(0)]}} \right) \right].\nonumber
\end{eqnarray}
}\end{prop}
\begin{proof} 
From (\ref{eq:W_1}) for $k=0,1,\ldots$ we have:
\begin{eqnarray}
& & W_1(t+k\,Q|x_0)  =  {|\Lambda_3(t+k\,Q)|^{1/2}\over 2\,\pi\, [1-\gamma^2(t+k\,Q)]}
\exp \left\{ -{[S(t+k\,Q)-x_0\,\gamma (t+k\,Q)]^2\over 2\,[1-\gamma^2(t+k\,Q)]}\right\} 
\label{eq:(4.28)}\\
& &\hspace*{1.2cm} \times  \left[ \exp \left\{ -{\sigma^2(t+k\,Q\,|\,x_0)
\over 2}\right\} 
-\sqrt{{\pi\over 2}}\,\sigma(t+k\,Q\,|\,x_0)\,\hbox{Erfc} \left( 
{\sigma(t+k\,Q\,|\,x_0)\over \sqrt{2}} \right) \right],\nonumber
\end{eqnarray}
where $\sigma(t|x_0)$ and $|\Lambda_3(t)|$ are defined in 
(\ref{eq:elementW_1}).  By virtue of (\ref{eq:Z(t)}), 
from (\ref{eq:asymperiodicboundary})  it follows:
\begin{eqnarray}
\lim_{k \to \infty} S(t+k\,Q) = S_0 + Z(t),\qquad 
\lim_{k \to \infty} \dot{S}(t+k\,Q) = \dot{Z}(t).
\label{eq:(4.29)}
\end{eqnarray}
Furthermore, due to (\ref{eq:asymconstantgamma1}) and (\ref{eq:Z(t)}), 
from (\ref{eq:elementW_1}) one has:
\begin{eqnarray}
\lim_{k \to \infty} |\Lambda_3(t+k\,Q\,|\,x_0)| = -\ddot{\gamma}(0),\qquad 
\lim_{k \to \infty} \sigma(t+k\,Q\,|\,x_0) = {\dot{Z}(t)\over \sqrt{
-\ddot{\gamma}(0)}}.
\label{eq:(4.30)}
\end{eqnarray}
Taking the limit as $k\to +\infty$ in (\ref{eq:(4.28)}), and making use of 
(\ref{eq:(4.29)}) and (\ref{eq:(4.30)}), we are then led to (\ref{eq:(4.27)}),  
which completes the proof. \end{proof}
\noindent
\begin{rem}\label{remark1}{\em
For all $t>0$ the function $R[Z(t)]$ defined in (\ref{eq:(4.27)}) is a positive,  
periodic function with period $Q$. Furthermore, $R[Z(t)]$ can also be written  as
\begin{eqnarray}
R[Z(t)] =  {1\over 2\, \pi\, \sqrt{[-\ddot{\gamma}(0)]}} 
\exp \left\{ - {[S_0+Z(t)]^2\over 2} \right\} \int_{\dot{Z}(t)}^{+\infty} \, 
[\xi - \dot{Z}(t)] \, \exp \left\{ - {\xi^2\over 2\,[- \ddot{\gamma}(0)]} 
\right\}\;d\xi.\nonumber\\ 
&&\label{eq:(4.31)}
\end{eqnarray}
}\end{rem}
\begin{proof}
Since $Z(t)$ is a periodic function of period $Q$, for all $k=0,1, \ldots$ one has 
$Z(t+k\,Q)=Z(t)$ and $\dot{Z}(t+k\,Q)=\dot{Z}(t)$. Hence,  from (\ref{eq:(4.27)}) 
it follows that $R[Z(t+k\,Q)]=R[Z(t)]$ $\;(k=0,1,\ldots)$, i.e. $R[Z(t)]$ is a 
periodic function with period $Q$.
\par 
We  shall now prove  that $R[Z(t)]>0$ for all $t>0$. To show it, let us re-write  
$R[Z(t)]$ as 
\begin{eqnarray}
R[Z(t)]= {\sqrt{-\ddot{\gamma}(0)}\over 2\,\pi} 
\exp \left\{ - {[S_0+Z(t)]^2\over 2} \right\}\,
A\biggl[{ \dot{Z}(t)\over\sqrt{2\,[-\ddot{\gamma}(0)]}}\biggr],
\label{eq:(4.32)}
\end{eqnarray}
with
\begin{eqnarray}
A(y):=e^{-y^2} - \sqrt{\pi}\;y\;{\rm Erfc}(y)
=e^{-y^2} - 2\;y\;\int_{y}^{+\infty} e^{-x^2} dx.
\label{eq:(4.33)}
\end{eqnarray}
We note that $A(y)>0$ for all $y\leq 0$. Furthermore, if $y>0$ the 
following inequality  (cf. \cite{Abramowitz72}) holds:
$$
{ e^{-y^2}\over y+\sqrt{y^2+2}}<\int_{y}^{+\infty} e^{-x^2} dx\leq 
{e^{-y^2}\over y+\sqrt{y^2+{\ds{4\over \pi}}} }\cdot
$$
Hence, one has
$$
e^{-y^2}\,{\sqrt{y^2+{\ds{4\over \pi}}}-y\over y+\sqrt{y^2+{\ds{4\over \pi}}}}
\leq A(y)<e^{-y^2}\,{\sqrt{y^2+2}-y\over y+\sqrt{y^2+2}}\qquad (y>0),
$$
which again implies $A(y)>0$. Therefore, for all $y\in{\bf R}$ there holds 
$A(y)>0$. Recalling (\ref{eq:(4.32)}), it  immediately follows  $R[Z(t)]>0$ for all 
$t>0$.
\par
Finally, to prove (\ref{eq:(4.31)}) we  note that 
\begin{eqnarray}
\int_{\dot{Z}(t)}^{+\infty}\left[\xi-\dot{Z}(t)\right]
\exp \left( - {\xi^2\over 2\,[- \ddot{\gamma}(0)]} \right)\;d\xi 
=[-\ddot{\gamma}(0)]\;A\biggl[{ \dot{Z}(t)\over\sqrt{2\,[-\ddot{\gamma}(0)}]}\biggr],
\label{eq:(4.34)}
\end{eqnarray}
with $A(y)$ defined in (\ref{eq:(4.33)}). Hence, making use (\ref{eq:(4.34)}) in 
(\ref{eq:(4.32)}), equation (\ref{eq:(4.31)}) immediately follows, 
as it had to be proved. \end{proof}
Since $Z(t)$ does not depend on $S_0$, from (\ref{eq:(4.27)}) it follows:
\begin{eqnarray}
\lim_{S_0\to +\infty}R\bigl[Z(t)\bigr]=0.
\label{eq:(4.35)}
\end{eqnarray}
%
\begin{prop}{\em
Let
\begin{eqnarray}
\alpha\equiv \alpha(S_0):= {1\over Q} \int_0^Q R[Z(\tau)]\;d\tau,
\label{eq:(4.36)}
\end{eqnarray}
with $R[Z(t)]$ defined in (\ref{eq:(4.27)}). Then, there exists a non-negative 
monotonically increasing function $\varphi(t)$ which is a solution of 
\begin{eqnarray}
\int_0^{\varphi(t)} R[Z(\tau)]\;d\tau = \alpha\,t,\qquad\qquad \forall t>0
\label{eq:(4.37)}
\end{eqnarray}
such that
\begin{eqnarray}
& &\varphi(0)=0,\nonumber\\
& &\lim_{t\to +\infty}\varphi(t)=+\infty,\label{eq:(4.38)}\\
& &\varphi(t+k\,Q) = \varphi(t) + k\,Q \qquad (k=0,1, \ldots).\nonumber
\end{eqnarray}
}\end{prop}
\begin{proof}
In Remark \ref{remark1} we have proved that $R[Z(t)]> 0,\;\forall t>0$; hence, from 
(\ref{eq:(4.36)}) it follows $\alpha>0$, and from (\ref{eq:(4.37)}) one has  
$\varphi(0)=0$ and $\varphi(t)> 0,\;\forall t>0$. Let $h(t)$ be any primitive 
function of $R[Z(t)].$ From (\ref{eq:(4.37)}) we have $h[\varphi(t)] = h(0) + \alpha\,t$.  
Since $R[Z(t)]> 0, \, \forall t>0$, $h(t)$ possesses an inverse, and hence
$\varphi(t)=h^{-1}[h(0)+\alpha\,t]$. Furthermore, since $\alpha > 0$, from 
(\ref{eq:(4.37)}) one has 
\begin{eqnarray}
{d\over dt}\,\varphi(t) = {\alpha\over R\left[Z(\varphi(t))\right]} > 0\qquad(t>0).
\label{eq:(4.39)}
\end{eqnarray}
Therefore, $\varphi(t)$ is a monotonically increasing function for all $t>0$. Furthermore, 
since $R[Z(t)]$ is a positive function, the second of  (\ref{eq:(4.38)}) 
holds. We now remark that from (\ref{eq:(4.37)}) one has:
\begin{eqnarray}
\int_0^{\varphi(t + k\,Q)}R[Z(\tau)]\;d\tau=
\int_0^{\varphi(t)} R[Z(\tau)]\;d\tau + \int_{\varphi(t)}^{\varphi(t + k\,Q)} 
R[Z(\tau)]\;d\tau = \alpha\,(t+k\,Q),
\label{eq:(4.40)}
\end{eqnarray}
or, due to (\ref{eq:(4.36)}),
\begin{eqnarray}
\int_{\varphi(t)}^{\varphi(t + k\,Q)} R[Z(\tau)]\;d\tau = 
k\,\alpha\, Q=k \int_0^{Q} R[Z(\tau)]\; d\tau = \int_0^{k\,Q} R[Z(\tau)]\;d\tau,
\label{eq:(4.41)}
\end{eqnarray}
where the last equality follows since $R[Z(t)]$ is a periodic function 
with period $Q$. Relation (\ref{eq:(4.41)}) finally implies the last of 
(\ref{eq:(4.38)}).  The proof is now complete.\end{proof}
\noindent
\begin{prop}\label{lemmaphi}{\em
For all $t>0$ one has
\begin{eqnarray}
& & {\it (i)}\quad\varphi\Bigl({t\over\alpha}\Bigr)>0,\nonumber\\
& & {\it (ii)}\quad
{d\over dt}\,\varphi \left( {t\over\alpha }\right)
={1\over R\biggl[Z\biggl(\varphi \biggl({\ds{t\over\alpha }}\biggr)\biggr)\biggr]}\,,
\nonumber\\
& & {\it (iii)}\quad\lim_{S_0\to +\infty}\,
\varphi\Bigl({t\over\alpha}\Bigr)=+\infty,\nonumber\\
& &{\it (iv)}\quad\lim_{S_0\to +\infty}\,\Bigl[\varphi\Bigl({t\over\alpha}\Bigr)
-\varphi\Bigl({\tau\over\alpha}\Bigr)\Bigr]=+\infty\qquad(0<\tau<t).\nonumber
\end{eqnarray}
}\end{prop}
\begin{proof}
Since $\varphi(t)$ is a non-negative function and $\alpha>0$,  condition {\it (i)} 
follows from (\ref{eq:(4.36)}), while from  (\ref{eq:(4.39)}) immediately one 
obtains {\it (ii)}. Making use of (\ref{eq:(4.35)}), from 
(\ref{eq:(4.36)}) we have
\begin{eqnarray}
\lim_{S_0\to +\infty}\alpha={1\over Q} \lim_{S_0\to +\infty}\int_0^Q R[Z(\tau)]\;d\tau
=0,
\label{eq:(4.43)}
\end{eqnarray}
that, due to the second of (\ref{eq:(4.38)}), implies {\it (iii)}. 
Finally, making use of the mean theorem of Calculus one has:
\begin{eqnarray}
\varphi\Bigl({t\over\alpha}\Bigr)-\varphi\Bigl({\tau\over\alpha}\Bigr)
={t-\tau\over\alpha}\;\dot{\varphi}\biggl({\xi\over\alpha}\biggr)
={t-\tau\over R\biggl[Z\biggl(\varphi\biggl({\ds{\xi\over\alpha}}
\biggr)\biggr)\biggr]}\qquad (0<\tau\leq\xi\leq t),
\label{eq:(4.43bis)}
\end{eqnarray}
where the last identity follows from {\it (ii)}. We note that, 
due to  (\ref{eq:(4.35)}), there holds:
\begin{eqnarray}
\lim_{S_0\to +\infty}R\biggl[Z\biggl(\varphi\biggl({t\over\alpha}\biggr)\biggr)
\biggr]=0,
\label{eq:(4.44)}
\end{eqnarray}
so that {\it (iv)} follows after taking the limit as $S_0\to +\infty$ 
in (\ref{eq:(4.43bis)}). The proof is thus complete.\end{proof}
The following theorem then holds.
%
%
\begin{thm}\label{thmperiodic}{\em
Let $S(t),t\geq 0$ be given in (\ref{eq:asymperiodicboundary})  
with $\varrho(t) \in C^1[0,+\infty)$ a bounded function such that 
(\ref{eq:Z(t)}) and (\ref{eq:intZ(t)}) hold. If the covariance function $\gamma(t)$ 
satisfies (\ref{eq:asymconstantgamma2})  and if 
\begin{eqnarray}
\lim_{S_0\to +\infty }\;
{\varrho \biggl(\varphi \biggl({\ds{t\over \alpha }}\biggr) \biggr) 
\over S_0+Z\biggl( \varphi \biggl( {\ds{t\over \alpha }}\biggr) \biggr) }=0,
\label{eq:(4.45)}
\end{eqnarray}
with $\varphi(t)$ defined in (\ref{eq:(4.37)}), then 
\begin{eqnarray}
\lim_{S_0\to +\infty }\left[ {d\over dt}\varphi \left( {t\over\alpha }
\right) \right] \,
g\left[ \varphi \left( {t\over\alpha }\right) \,\,\Big|\,\,x_0\right]
=e^{-t},\qquad x_0<S(0).  
\label{eq:(4.46)}
\end{eqnarray}
}\end{thm}
\begin{proof}
From {\it (i)} and {\it (iii)} of Proposition \ref{lemmaphi} it follows that 
$\varphi\bigl(t/\alpha\bigr)$ can 
be viewed as a scaled time. Changing $t$ to  
$\varphi\bigl(t/\alpha \bigr)$ in (\ref{eq:FPTGaussiandensity}),  
we then obtain:
\begin{eqnarray*}
&&\hspace*{-0.5cm}\biggl[ {d\over dt}\,\varphi \biggl( {t\over \alpha}
\biggr) \biggr] \,
g\biggl[\varphi \biggl( {t\over \alpha }\biggr) \;\Big|\; x_0\biggr] 
=\biggl[ {d\over dt}\,\varphi \biggl( {t\over\alpha }\biggr) \biggr] 
W_1\biggl[ \varphi \biggl( {t\over\alpha }\biggr) \;\Big|\; x_0\biggr] 
+\biggl[ {d\over dt}\varphi \biggl( {t\over\alpha }\biggr) \biggr]\\
& & \times\sum_{i=1}^{\infty }(-1)^i
\int_0^{\varphi \bigl( {t\over\alpha }\bigr) }dt_1
\int_{t_1}^{\varphi \bigl( {t\over \alpha }\bigr)}dt_2
\cdots \int_{t_{i-1}}^{\varphi \bigl( {t\over \alpha }\bigr)}
W_{i+1}\biggl( t_1,\cdots ,t_i,{\varphi \biggl( {t\over \alpha }
\biggr) }\;\Big|\; x_0\biggr)\; dt_i\,.
\end{eqnarray*}
Hence: 
\begin{eqnarray}
& &\biggl[ {d\over dt}\,\varphi \biggl( {t\over \alpha }\biggr) \biggr] \,
g\biggl[\varphi \biggl( {t\over \alpha }\biggr) \;\Big|\; x_0\biggr] 
=\biggl[ {d\over dt}\,\varphi \biggl( {t\over\alpha }\biggr) \biggr] 
W_1\biggl[ \varphi \biggl( {t\over\alpha }\biggr) \;\Big|\; x_0\biggr] 
\label{eq:(4.49)}\\
& &\hspace*{1cm}+\biggl[ {d\over dt}\varphi \biggl( {t\over \alpha }\biggr) \biggr] 
\;\sum_{i=1}^{\infty }(-1)^i\int_0^{t}\biggl[ {d\over d\tau_1}
\varphi \biggl( {\tau_1\over\alpha }\biggr) \biggr] d\tau_1
\int_{\tau_1}^{t}\biggl[ {d\over d\tau_2}
\varphi \biggl( {\tau_2\over\alpha }\biggr) \biggr] d\tau_2\cdots  
\nonumber\\
&&\hspace*{1cm}\cdots\int_{\tau_{i-1}}^{t}\biggl[ {d\over d\tau_i}
\varphi \biggl( {\tau_i\over\alpha }\biggr) \biggr] 
W_{i+1}\biggl( \varphi \biggl( {\tau_1\over\alpha }\biggr) ,\ldots,
\varphi \biggl( {\tau_i\over\alpha }\biggr) ,
{\varphi \biggl( {t\over\alpha}\biggr) }\;\Big|\;x_0\biggr)\;d\tau_i\,,\nonumber 
\end{eqnarray}
after having performed the change of variables 
$\tau_i=\varphi \left( {\ds{t_i\over\alpha }}\right)$.  
Due to {\it (ii)} of Proposition \ref{lemmaphi}, (\ref{eq:(4.49)}) can also be  
written as
\begin{eqnarray}
& &\hspace*{-1.5cm}\biggl[ {d\over dt}\,\varphi \biggl( {t\over \alpha}
\biggr) \biggr] \,
g\biggl[\varphi \biggl( {t\over \alpha }\biggr) \;\Bigl|\; x_0\biggr] 
={W_1\biggl[ \varphi \biggl( {\ds{t\over\alpha }}\biggr) \;\Bigl|\; x_0\biggr]
\over  R\biggl[Z\biggl(\varphi \biggl({\ds{t\over\alpha }}\biggr)\biggr)\biggr]}
+ \sum_{i=1}^{\infty }(-1)^i\int_0^{t} d\tau_1\int_{\tau_1}^t d\tau_2\cdots
\label{eq:(4.49bis)}\\
& &\hspace*{0.5cm}\cdots  \int_{\tau_{i-1}}^t \,
{W_{i+1}\biggl( \varphi \biggl( {\ds{\tau_1\over\alpha }}\biggr) ,\ldots,
\varphi \biggl( {\ds{\tau_i\over\alpha }}\biggr) ,
{\varphi \biggl( {\ds{t\over\alpha}}\biggr) }\;\Bigl|\;x_0\biggr)\over
R\biggl[Z\biggl(\varphi \biggl({\ds{\tau_1\over\alpha }}\biggr)\biggr)\biggr]
\cdots R\biggl[Z\biggl(\varphi \biggl({\ds{\tau_i\over\alpha }}\biggr)\biggr)\biggr]\,
R\biggl[Z\biggl(\varphi \biggl({\ds{t\over\alpha }}\biggr)\biggr)\biggr]
}\;d\tau_i\, .\nonumber 
\end{eqnarray}
Let us now prove that for  $n=1,2,\ldots$ there holds: 
\begin{eqnarray}
\hspace*{-0.5cm}\lim_{S_0\to +\infty }
{ W_n\biggl( \varphi \biggl( {\ds{\tau_1\over\alpha }}\biggr) ,\ldots,
\varphi \biggl( {\ds{\tau_n\over\alpha} }\biggr) \;\Big|\;x_0\biggr)
\over R\biggl[Z\biggl(\varphi \biggl({\ds{\tau_1\over\alpha }}\biggr)
\biggr)\biggr]\cdots
R\biggl[Z\biggl(\varphi \biggl({\ds{\tau_n\over\alpha }}\biggr)\biggr)\biggr] }=1
\qquad(0<\tau_1<\tau_2<\ldots <\tau_n). 
\label{eq:(4.50)}
\end{eqnarray}
For simplicity of notation, we set:
\begin{eqnarray}
\vartheta_1=\varphi \biggl({\tau_1\over\alpha }\biggr),\ldots, 
\vartheta_n=\varphi \biggl({\tau_n\over\alpha }\biggr)\qquad
(0 < \vartheta_1 < \vartheta_2 < \ldots < \vartheta_n).
\label{eq:(4.50ter}
\end{eqnarray}
From (\ref{eq:(4.31)}), for $n=1,2,\ldots$ it then follows
\begin{eqnarray}
& &\hspace*{-0.5cm}R\Bigl[Z\Bigl(\vartheta_1\Bigr)\Bigr]
\cdots R\biggl[Z\Bigl(\vartheta_n\Bigr)\Bigr]
={1\over\left( 2\,\pi \right)^n\left[ -\ddot{\gamma}(0)
\right]^{n/2}}\,
\exp \biggl\{ -{1\over 2}\sum\limits_{i=1}^n\left[ S_0+Z\left( 
\vartheta_i \right) \right]^2\biggr\} 
\label{eq:(4.52)}\\
& &\hspace*{0.3cm}\times\int_{\dot{Z}\left( \vartheta_1 \right) }^{+\infty }d\xi_1
\int_{\dot{Z}\left(\vartheta_2 \right)}^{+\infty }d\xi_2
\cdots\int_{\dot{Z}\left( \vartheta_n \right)}^{+\infty}
\prod_{i=1}^n\left[ \xi_i-\dot{Z}\left( 
\vartheta_i \right) \right] 
\exp \biggl\{ -{1\over 2\,[-\ddot{\gamma}(0)]}\sum_{i=1}^n\xi_i^2\biggr\}\;d\xi_n\,. 
\nonumber
\end{eqnarray}
Hence,
\begin{eqnarray}
&&{W_n\left( \vartheta_1 ,\cdots,\vartheta_n \;\big|\;x_0\right)
\over R\left[Z\left(\vartheta_1\right)\right]
\cdots R\left[Z\left(\vartheta_n\right)\right]}
=\left( 2\,\pi \right)^n\left[ -\ddot{\gamma}(0)\right]^{n/2}
\label{eq:(4.52)bis}\\
&&\times\,{ \exp \biggl\{ {\ds{1\over2}}\sum\limits_{i=1}^n
\bigl[ S_0+Z\bigl( \vartheta_i \bigr) \bigr]^2\biggr\}\,
W_n\left( \vartheta_1 ,\cdots,\vartheta_n \;\big|\;x_0\right)\over
{\ds\int_{\dot{Z}\left( \vartheta_1 \right) }^{+\infty }}d\xi_1
{\ds\int_{\dot{Z}\left( \vartheta_2 \right) }^{+\infty }}d\xi_2
\cdots {\ds\int_{\dot{Z}\left( \vartheta_n \right)}^{+\infty}}
\prod\limits_{i=1}^n\left[ \xi_i-\dot{Z}\left( 
\vartheta_i \right) \right] \exp \biggl\{ -{\ds{1\over 2\,[-\ddot{\gamma}(0)]}}
\,\sum\limits_{i=1}^n\xi_i^2\biggr\}\; d\xi_n}\nonumber
\end{eqnarray}
where, due to (\ref{eq:def2W_n}),  one has
\begin{eqnarray}
& &\hspace*{-1.0cm}W_n\left( \vartheta_1 ,
\cdots,\vartheta_n \;\big|\;x_0\right)
={ 1\over (2\,\pi)^n \left| \Lambda_{2n+1} \left( \vartheta_1,\ldots,
\vartheta_n \right) \right|^{1/2} } 
\int_{\dot{\psi}\left(\vartheta_1\right)}^{+\infty} d\xi_1
\int_{\dot{\psi}\left(\vartheta_2\right)}^{+\infty} d\xi_2
\cdots 
\label{eq:(4.53)}\\ 
& &\hspace{-0.5cm}\cdots\int_{\dot{\psi}\left(\vartheta_n\right)}^{+\infty} 
d\xi_n\prod_{i=1}^n \left[ \xi_i-\dot{\psi}
\left( \vartheta_i\right) \right]\,
H \left( \vartheta_1,\ldots,\vartheta_n
;\xi_1,\ldots,\xi_n \right)\, K\left( \vartheta_1,
\ldots,\vartheta_n;\xi_1,\ldots,\xi_n \right)\nonumber 
\end{eqnarray}
with $\psi$, $H$ and $K$ defined in (\ref{eq:componentW_n}).
\par\noindent
We now note that (see Appendix II)
\begin{eqnarray}
& & \lim_{S_0 \to +\infty} \, \Big| \, \Lambda_{2n+1} \left(\vartheta_1, 
\ldots,\vartheta_n \right) \, \Big| \, = \bigl[-\ddot{\gamma}(0)\bigr]^n,  
\label{eq:appendix1}\\
& & \lim_{S_0 \to +\infty}\, 
H\left( \vartheta_1,\ldots,\vartheta_n;\xi_1,\ldots,\xi_n \right)
=\exp \biggl\{ - {1\over 2\,[-\ddot{\gamma}(0)]} \sum_{i=1}^n
\xi_i^2 \biggr\},\label{eq:appendix2}\\
&& \lim_{S_0 \to +\infty} \biggl[\exp \biggl\{ {1\over2}\sum\limits_{i=1}^n
\bigl[ S_0+Z\bigl( \vartheta_i \bigr) \bigr]^2\biggr\}\,
K\left( \vartheta_1,\ldots,\vartheta_n;\xi_1,\ldots,\xi_n \right)\biggr]=1.\qquad
\label{eq:(4.54)quat}
\end{eqnarray}
Taking the limit as $S_0 \to +\infty$ in (\ref{eq:(4.52)bis}) and 
recalling (\ref{eq:appendix1}), (\ref{eq:appendix2}) and  
(\ref{eq:(4.54)quat}),  for all $n=1,2,\ldots$  one has
\begin{eqnarray}
&&\hspace*{-0.7cm}\lim_{S_0\to +\infty }
{W_n\left( \vartheta_1 ,\cdots,\vartheta_n \;\big|\;x_0\right)
\over R\left[Z\left(\vartheta_1\right)\right]
\cdots R\left[Z\left(\vartheta_n\right)\right]}
\label{eq:(4.55)}\\
&&\hspace*{-0.5cm}=\lim_{S_0\to +\infty }{ 
{\ds\int_{\dot{\psi}\left(\vartheta_1\right)}^{+\infty}} d\xi_1
{\ds\int_{\dot{\psi}\left(\vartheta_2\right)}^{+\infty}} d\xi_2
\cdots {\ds\int_{\dot{\psi}\left(\vartheta_n\right)}^{+\infty}} 
\prod\limits_{i=1}^n \Bigl[ \xi_i-\dot{\psi}
\left( \vartheta_i\right) \Bigr]\,\exp \biggl\{ - {\ds{1\over 2\,[-\ddot{\gamma}(0)]}} 
\sum\limits_{i=1}^n\xi_i^2 \biggr\}\;d\xi_n\over
{\ds\int_{\dot{Z}\left( \vartheta_1 \right) }^{+\infty }}d\xi_1
{\ds\int_{\dot{Z}\left( \vartheta_2 \right) }^{+\infty }}d\xi_2
\cdots {\ds\int_{\dot{Z}\left( \vartheta_n \right)}^{+\infty}}
\prod\limits_{i=1}^n\left[ \xi_i-\dot{Z}\left( 
\vartheta_i \right) \right] \exp \left\{ -{\ds{1\over 2\,[-\ddot{\gamma}(0)]}}
\,\sum\limits_{i=1}^n\xi_i^2\right\}\;d\xi_n }\nonumber\\
&& \hspace*{-0.5cm}=\lim_{S_0\to +\infty }\prod_{i=1}^n
{U\left( \vartheta_i\;\big|\;x_0\right)\over V\left( \vartheta_i \right)}
\qquad\qquad \qquad(0<\vartheta_1<\vartheta_2<\ldots <\vartheta_n),
\nonumber 
\end{eqnarray}
where we have set:
\begin{eqnarray}
& & U\left( \vartheta_i\;\vline\;x_0 \right)=
\int_{\dot{\psi}\left(\vartheta_i\right)}^{+\infty} 
\left[ \xi-\dot{\psi}\left( \vartheta_i\right) \right]\,
\exp \left\{ - {\xi^2\over 2\,[-\ddot{\gamma}(0)]} \right\}\;d\xi\, ,\nonumber\\
& &\label{eq:(4.56)} \\
& & V\left( \vartheta_i\right)=
\int_{\dot{Z}\left( \vartheta_i\right) }^{+\infty } 
\left[ \xi -\dot{Z}\left( \vartheta_i \right) \right]\,
\exp \left\{ -{\xi^2\over 2\,[-\ddot{\gamma}(0)]}\right\}\;d\xi\,.\nonumber
\end{eqnarray}
We now note that from (\ref{eq:(4.33)}),  (\ref{eq:(4.34)}) and the first of  
(\ref{eq:componentW_n}) it follows
\begin{eqnarray}
& &\hspace*{-1.0cm} U\left( \vartheta_i\;\vline\;x_0 \right)=
-\ddot{\gamma}(0)\,\exp \biggl\{ - {\bigl[\dot{S}\left(
\vartheta_i\right)-x_0\,\dot{\gamma}\left(\vartheta_i\right)
\bigr]^2\over 2\,[-\ddot{\gamma}(0)]} \biggr\}-  { \sqrt{\pi\,[-\ddot{\gamma}(0)]\over 2}} 
\nonumber\\
& &\hspace*{1.0cm}\times\left[
\dot{S}\left(\vartheta_i\right)\, 
-x_0\,\dot{\gamma}\left(\vartheta_i\right)
\right]\,{\rm Erfc} \, \left({\dot{S}\left(
\vartheta_i\right)\, 
-x_0\,\dot{\gamma}\left(\vartheta_i\right)
\over \sqrt{2\,[-\ddot{\gamma}(0)] }} \right) ,\nonumber\\
& &\label{eq:(4.57)} \\  
& &\hspace*{-1.0cm} V\left( \vartheta_i\right)=
-\ddot{\gamma}(0)\,\exp \biggl\{ - {\bigl[\dot{Z}\left(
\vartheta_i\right)\bigr]^2\over 2\,[-\ddot{\gamma}(0)]} \biggr\}
-  { \sqrt{\pi\,[-\ddot{\gamma}(0)]\over 2}}  \, \dot{Z}\left(
\vartheta_i\right) \, {\rm Erfc} \, \left({\dot{Z}\left(\vartheta_i\right)
\over \sqrt{2\,[-\ddot{\gamma}(0)] }} \right).\nonumber
\end{eqnarray}
Furthermore, recalling (\ref{eq:asymconstantgamma2}) and (\ref{eq:(4.29)}), one has
$$
\lim_{S_0\to +\infty }
{U\left( \vartheta_i\;\vline\;x_0\right)
\over V\left( \vartheta_i \right)}=1 \qquad(i=1,2,\ldots,n), 
$$
so that (\ref{eq:(4.50)}) immediately follows  from (\ref{eq:(4.55)}). 
Due to (\ref{eq:(4.50)}), taking the limit as $S_0\to +\infty $ in 
(\ref{eq:(4.49)}), one obtains: 
\begin{eqnarray}
\lim_{S_0\to +\infty }\left[ {d\over dt}\varphi \left( {t\over \alpha }
\right) \right] g\left[ \varphi \left( {t\over \alpha }\right) \;
\Big|\;x_0\right]=1 + \sum_{i=1}^\infty (-1)^i 
\int_0^{t} d\tau_1 \int_{\tau_1}^{t} d\tau_2 \cdots \int_{\tau_{i-1}}^{t} d\tau_i. 
\label{eq:(4.58)}
\end{eqnarray}
Equation (\ref{eq:(4.46)}) finally follows 
from (\ref{eq:(4.58)}).  The proof is thus complete.\end{proof}
%
\begin{cor}{\em 
Under the assumption of Theorem \ref{thmperiodic}, for 
$S_0\to +\infty$ one has: 
\begin{eqnarray}
g(t \,|\, x_0) \sim R[Z(t)]\,\exp\biggl\{ -\int_0^t R[Z(\tau)]\;d\tau\biggr\},
\qquad\forall t>0, 
\label{eq:(4.59)}
\end{eqnarray}
with $R[Z(t)]$ defined in (\ref{eq:(4.27)}).
}\end{cor}
\begin{proof}
From {\it (ii)} of Proposition \ref{lemmaphi} and from (\ref{eq:(4.46)}) we obtain:
\begin{eqnarray}
\lim_{S_0\to +\infty }
{g\biggl[ \varphi \biggl( {\ds{t\over\alpha }}\biggr) \,\,\Big|\,\,x_0\biggr]
\over R\biggl[Z\biggl(\varphi \biggl({\ds{t\over\alpha }}\biggr)\biggr)\biggr]}
=e^{-t},\qquad x_0<S(0).  
\label{eq:(4.60)}
\end{eqnarray}
Hence, using (\ref{eq:(4.37)}), we have asymptotically:
$$
g\biggl[ \varphi \biggl( {\ds{t\over\alpha }}\biggr) \,\,\Big|\,\,x_0\biggr]
\sim R\biggl[Z\biggl(\varphi \biggl({\ds{t\over\alpha }}\biggr)\biggr)\biggr]\,e^{-t}
=R\biggl[Z\biggl(\varphi \biggl({\ds{t\over\alpha }}\biggr)\biggr)\biggr]\,
\exp\biggl\{-\int_0^{\varphi(t/\alpha)} R[Z(\tau)]\;d\tau\biggr\}.
$$
The asymptotic formula (\ref{eq:(4.59)}) then follows. The asymptotic non-homogeneous 
exponential behavior of the FPT density has thus been proved.
\end{proof}
Note that (\ref{eq:(4.59)}) can also be written as
\begin{eqnarray}
g(t \,|\, x_0) \sim \beta(t)\,e^{-\alpha\,t},
\label{eq:(4.60bis)}
\end{eqnarray}
where $\beta(t)$ is a periodic function of period $Q$ given by
\begin{eqnarray}
\beta(t)=R[Z(t)]\,\exp\biggl\{\alpha\,t- \int_0^t R[Z(\tau)]\;d\tau\biggr\},
\label{eq:(4.60ter)}
\end{eqnarray}
with $\alpha$ defined in (\ref{eq:(4.36)}). Indeed, since $R[Z(t)]$ is a periodic 
function of period $Q$, due to (\ref{eq:(4.36)}) and (\ref{eq:(4.37)}), one has:
\begin{eqnarray}
&&\beta(t+nQ)=R[Z(t+nQ)]\,\exp\biggl\{\alpha\,(t+nQ)- \int_0^{t+nQ} 
R[Z(\tau)]\;d\tau\biggr\}\label{eq:(4.61)}\\
&&\hspace*{0.5cm}=\beta(t)\,\exp\biggl\{\alpha\,nQ- \int_t^{t+nQ} 
R[Z(\tau)]\;d\tau\biggr\}\nonumber\\
&&\hspace*{0.5cm}=\beta(t)\,\exp\biggl\{\int_0^{\varphi(nQ)} R[Z(\tau)]\;d\tau
- \int_0^{nQ} R[Z(\tau)]\;d\tau\biggr\}=\beta(t).\nonumber
\end{eqnarray}
\par
Notice that in (\ref{eq:(4.60bis)}) is an asymptotic expression of $g(t \,|\, x_0)$ 
of the same form as (\ref{eq:(1.1)}).
%
\begin{figure}[thb]
 \epsfxsize=11cm
\centerline{\epsfbox{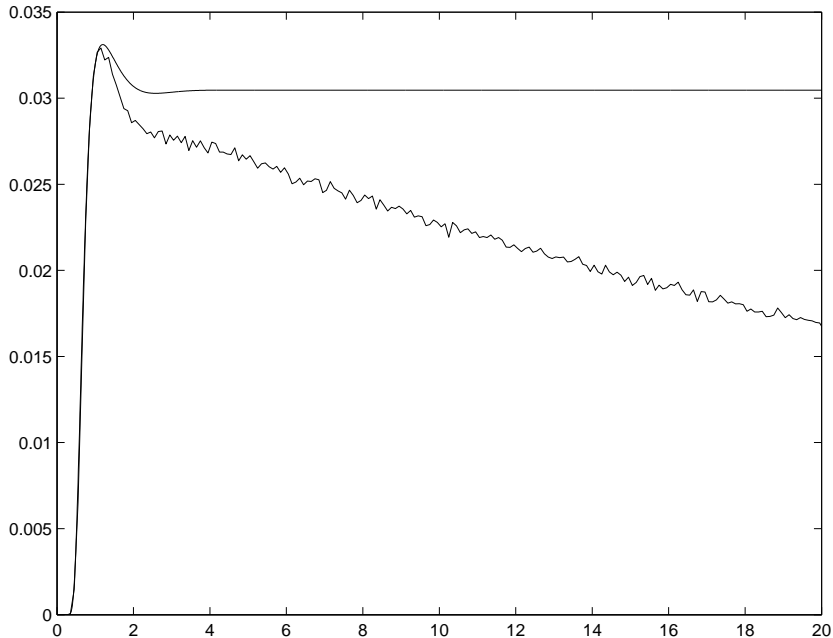}}
\caption {Plots of the function $W_1(t|0)$ and of the simulated FPT 
density $\tilde{g}(t)$ for a zero-mean stationary Gaussian process 
originating at $x_0=0$ having covariance (\ref{es:gencovariance}) 
with $a=\omega=1$ and for the constant boundary $S(t)=2$.}
\label{figure:g_W1_2}
\end{figure}
%
\begin{figure}[thb]
\epsfxsize=10cm 
\centerline{\epsfbox{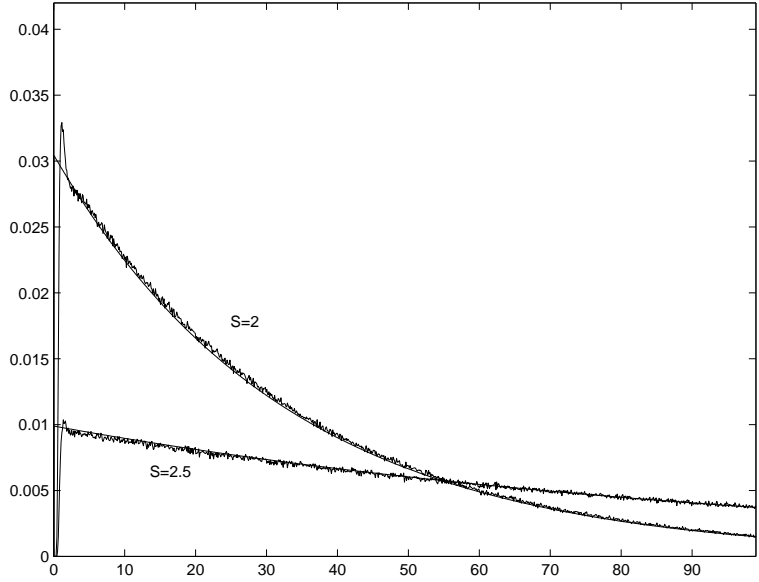}}
\caption{For the  boundaries $S(t)\equiv S_0=2$ and $S(t)\equiv S_0=2.5$, the asymptotic 
exponential approximation $R(S_0)\,\exp\{-R(S_0)\,t\}$ and the simulated FPT 
density $\tilde{g}(t)$  are plotted for the Gaussian process of 
Figure \ref{figure:g_W1_2}.}
\label{figure:g_2-2.5}
\end{figure}
%
\begin{figure}[thb]
\epsfxsize=10cm
\centerline{\epsfbox{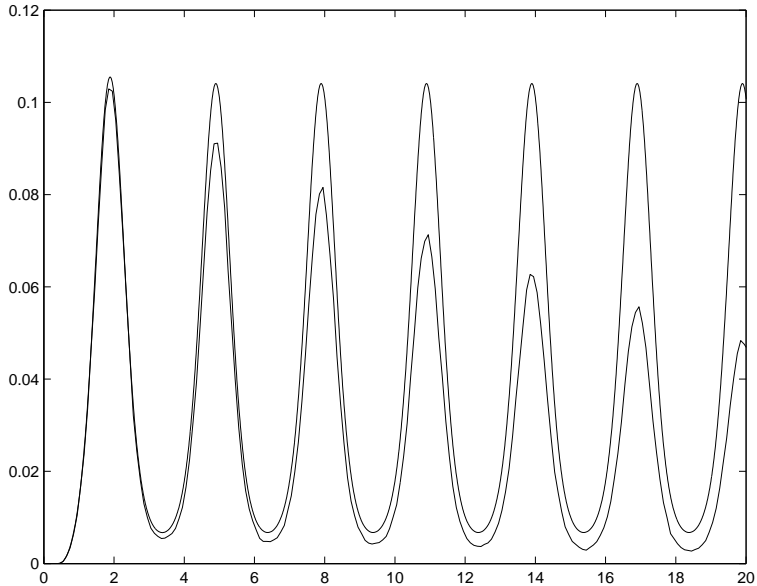}}
\caption{Same as in Figure 1 for the periodic boundary 
$S(t)=2+0.5\,\sin(2\,\pi\, t/3)$.}
\label{figure:g_W1_2_0.5_3}
\end{figure}
%
\begin{figure}[thb]
\epsfxsize=10cm \centerline{\epsfbox{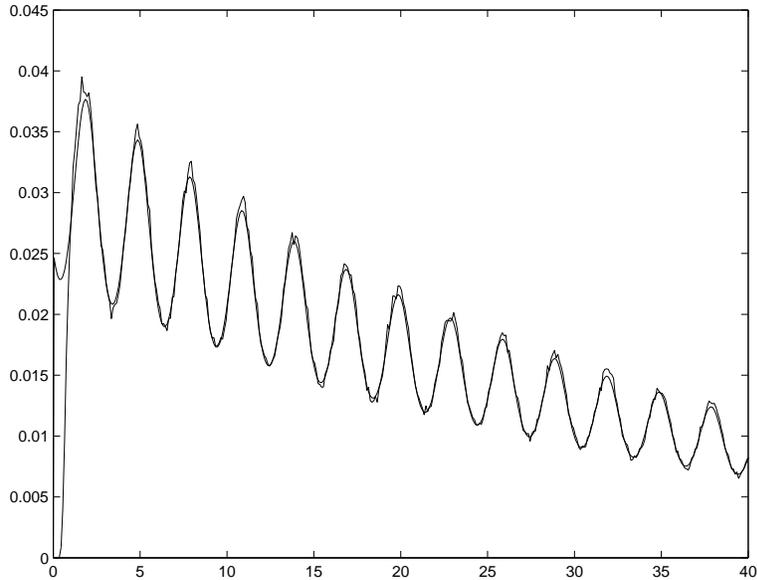}}
\caption{For the same process as in Figure \ref{figure:g_W1_2} and for 
$S(t)=2+0.1\,\sin(2\,\pi\,\,t/3)$,  the function 
$R[Z(t)]\,\exp\{-\int_0^tR[Z(\tau)]\,d\tau\}$ is 
plotted together with the simulated FPT density $\tilde{g}(t)$.}
\label{figure:g_2_0.1_3}
\end{figure}
%
\begin{figure}[thb]
\epsfxsize=10cm \centerline{\epsfbox{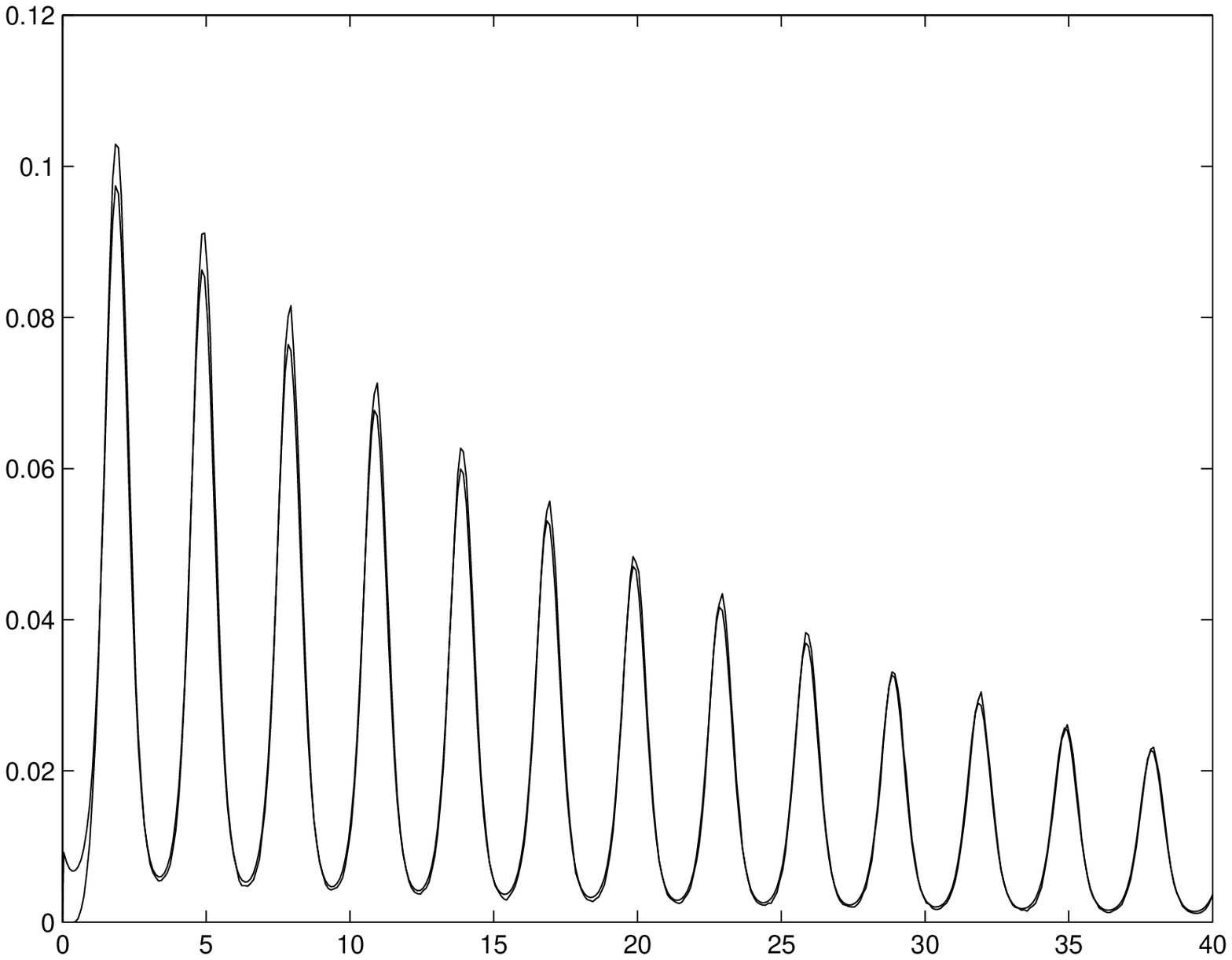}}
\caption{Same as in Figure \ref{figure:g_2_0.1_3} for the 
boundary $S(t)=2+0.5\,\sin(2\,\pi\,t/3)$.}
\label{figure:g_2_0.5_3}
\end{figure}
%
\begin{figure}[thb]
\epsfxsize=10cm \centerline{\epsfbox{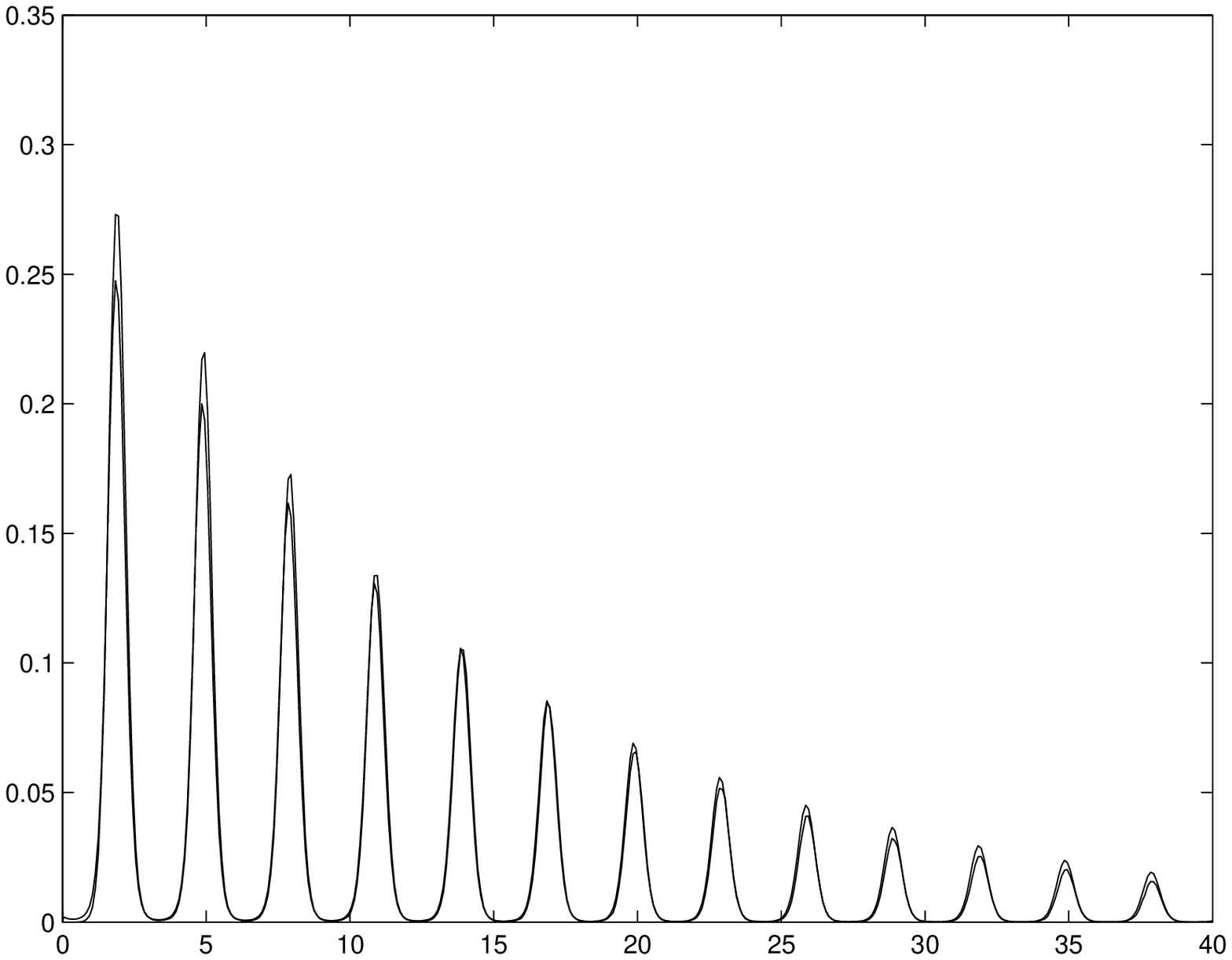}}
\caption{Same as in Figure \ref{figure:g_2_0.1_3} for the 
boundary $S(t)=2+\sin(2\,\pi\,t/3)$.}
\label{figure:g_2_1_3}
\end{figure}
%
\begin{figure}[thb]
\epsfxsize=10cm \centerline{\epsfbox{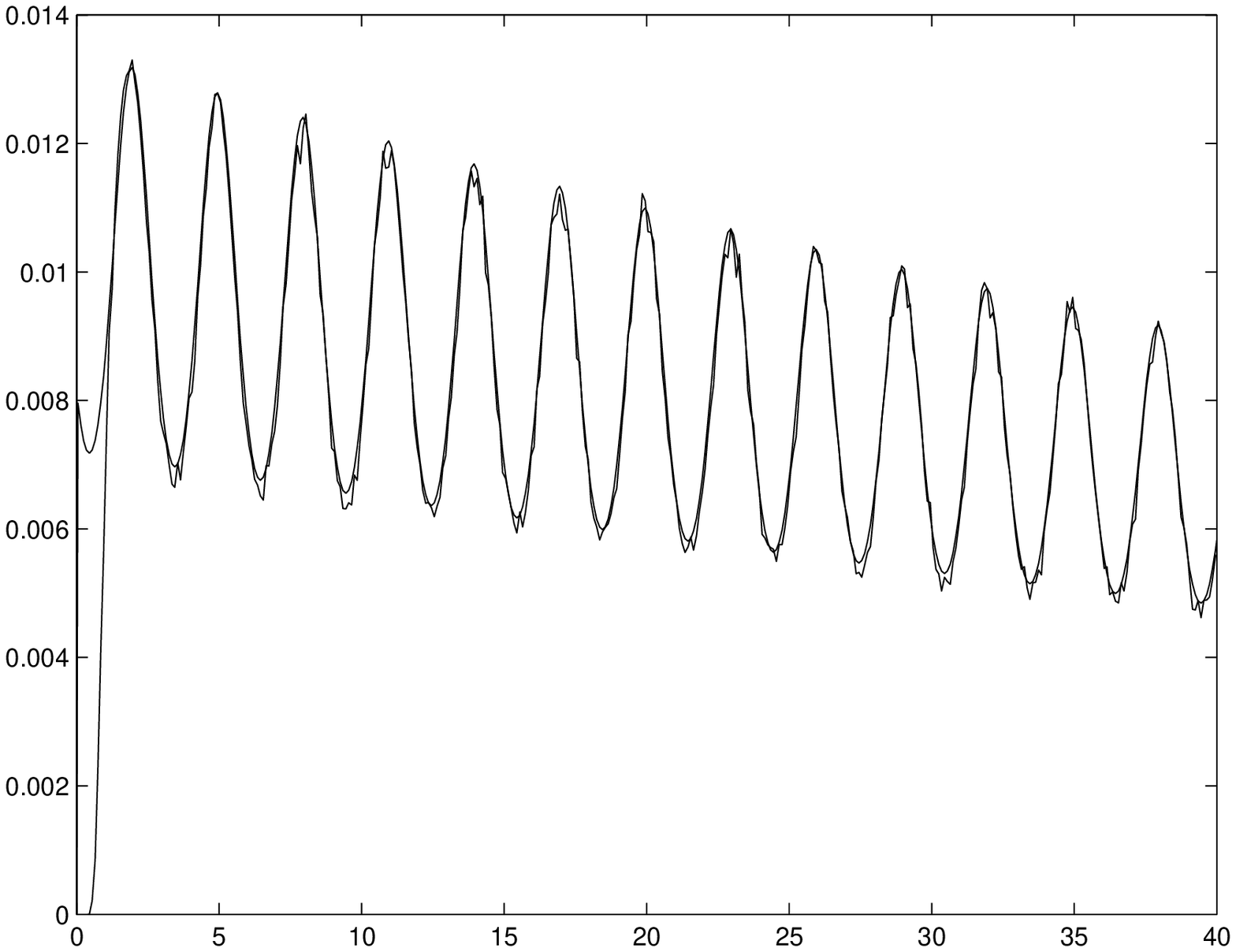}}
\caption{Same as in Figure \ref{figure:g_2_0.1_3} for the 
boundary $S(t)=2.5+0.1\,\sin(2\,\pi\,t/3)$.}
\label{figure:g_2.5_0.1_3}
\end{figure}
\par\noindent
\section{Analysis of a special case}
\setcounter{equation}{0}
The purpose of this Section is to analyze the behavior of the FPT 
pdf for a stationary Gaussian process of concrete interest for certain 
applications.
\par
Let $\left\{X(t), t\geq 0 \right\}$ be the stationary 
Gaussian process originating at $x_0=0$, with zero mean and damped 
oscillatory covariance \cite{Stratonovich63}: 
\begin{eqnarray}
\gamma (t)= e^{-a\,|t|} \,\Bigl[\cos\bigl(\omega\, t)+{a\over\omega}
\,\sin\bigl(\omega\cdot |t|)\Bigr]\qquad 
(t\in{\bf R}),
\label{es:gencovariance}
\end{eqnarray}
where $a$ and $\omega$ are positive real numbers. Functions of form 
(\ref{es:gencovariance})  can often be conveniently used to approximate 
experimental covariance functions that, starting from a unit initial maximum 
amplitude,  asymptotically tend to zero  with an exponential envelope.
From (\ref{es:gencovariance}) we see that $\gamma (0)=1$. Furthermore, 
$\dot{\gamma}(0)=0$ and $\ddot{\gamma}(0)=-(a^2+\omega^2)<0$
since for $t>0$ there holds:
\begin{eqnarray}
&&\dot{\gamma}(t)=-\dot{\gamma}(-t)=
-{e^{-a\,t}\over\omega}\,(a^2+\omega^2)\,\sin(\omega\, t),\nonumber\\
&&\label{es:dercovariance}\\
&&\ddot{\gamma}(t)=\ddot{\gamma}(-t)
={e^{-a\,t}\over\omega}\,(a^2+\omega^2)\,\Bigl[a\,\sin(\omega\, t)
-\omega\,\cos(\omega\, t)\Bigr].\nonumber
\end{eqnarray}
\par
We finally note that $\lim\limits_{t\to +\infty}\gamma(t)=0$, 
$\lim\limits_{t\to +\infty}\dot{\gamma}(t)=0$ 
and $\lim\limits_{t\to +\infty}\ddot{\gamma}(t)=0$. 
\subsection{Constant boundary}
For the constant boundary 
$S(t)=S_0$, from  (\ref{eq:R(S_0)}) one has
\begin{eqnarray}
R(S_0)={\sqrt{a^2+\omega^2}\over 2\,\pi}\; e^{\displaystyle{-S_0^2/2}}. 
\label{es:R(S_0)}
\end{eqnarray}
All forthcoming figures refer to the case $a=\omega=1$ in (\ref{es:gencovariance}). 
Figure \ref{figure:g_W1_2} shows the plots  of $W_1(t|x_0)$ given by (\ref{eq:W_1}) 
and of the simulated FPT density $\tilde{g}(t)$ as functions of $t$ for the constant 
boundary $S(t)=2$. The first-order approximation $W_1(t|x_0)$ is 
seen to provide an upper bound to the FPT pdf in (\ref{eq:FPTGaussiandensity}), 
being a good approximation of $g$ only for small values of $t$. 
From  (\ref{eq:R(S_0)}) it follows that, as $t$ increases,  $W_1(t|x_0)$  
approaches the constant value $R(2)=e^{-2}/(\pi\,\sqrt{2})=0.0304611$.
Making use of (\ref{eq:comp_g(t)}), we expect that $R(S_0)\,\exp\{-R(S_0)\,t\}$ 
provides a good approximation of the  simulated FPT density $\tilde{g}(t)$ 
as $S_0\to +\infty$. Indeed, Figure \ref{figure:g_2-2.5} shows the plots 
of the simulated FPT density $\tilde{g}(t)$ for the constant boundary 
$S(t)\equiv S_0=2$  and of the function $R(2)\,\exp\{-R(2)\,t\}$. 
Figure \ref{figure:g_2-2.5} also shows the plots of the simulated FPT 
density $\tilde{g}(t)$ for the constant 
boundary $S(t)\equiv S_0=2.5$ and the function $R(2.5)\,\exp\{-R(2.5)\,t\}$  
with $R(2.5)=e^{-3.125}/(\pi\,\sqrt{2})=0.00988928$. We note that  starting 
from rather small times,  $\tilde{g}(t)$ is susceptible of an excellent exponential 
approximation already for positive boundaries of the order of a couple of units. 
\subsection{Periodic boundary}
For the same stationary Gaussian process, we now assume that    
$S(t)=S_0+B\,\sin(2\,\pi\, t/Q)$. In this case 
$Z(t)\equiv\varrho(t)=B\,\sin(2\,\pi\, t/Q)$ and  condition  
(\ref{eq:intZ(t)}) is satisfied. From (\ref{eq:(4.27)}) one has:
\begin{eqnarray}
&&R[Z(t)]=
{\sqrt{a^2+\omega^2}\over 2\,\pi} \exp \biggl\{ - {S^2(t)\over 2} \biggr\}\, 
\biggl[ \exp \biggl\{ - {[\dot{Z}(t)]^2\over 2\,(a^2+\omega^2)} 
\biggr\}\nonumber \\ 
&&\hspace*{4cm}-  \sqrt{{\pi\over 2\,(a^2+\omega^2)}} \, \dot{Z}(t) \, {\rm Erfc} 
\, \biggl({\dot{Z}(t)\over\sqrt{2\,(a^2+\omega^2)}}  \biggr) \biggr].
\label{es:R[Z(t)]}
\end{eqnarray}
From Remark \ref{remark1} it follows that the function (\ref{es:R[Z(t)]})  is  
positive and periodic with period $Q$.
\par 
Figure \ref{figure:g_W1_2_0.5_3} shows the plots of  $W_1(t|0)$  
and of the simulated FPT density $\tilde{g}(t)$ for the 
periodic boundary $S(t)=2+0.5\,\sin(2\,\pi\, t/3)$. We note again that 
the first-order approximation $W_1$ provides an upper bound to the FPT pdf given by  
(\ref{eq:FPTGaussiandensity}), though being a good approximation of $g$ only for small 
values of $t$. Furthermore, as  $k$ increases   $W_1(t+k\,Q\,|\,0)$ 
tends to the function (\ref{es:R[Z(t)]}), i.e. as  $t$ increases 
$W_1(t|0)$ becomes periodic with the same period of the boundary.
Since 
$$
\lim_{S_0\to +\infty }\;
{\varrho \biggl(\varphi \biggl({\ds{t\over \alpha }}\biggr) \biggr) 
\over S_0+Z\biggl( \varphi \biggl( {\ds{t\over \alpha }}\biggr) \biggr) }
=\lim_{S_0\to +\infty }\;{B\,\sin\Bigl(2\,\pi\, {\ds{\varphi(t/\alpha)\over Q}}\Bigr)
\over S_0+B\,\sin\Bigl(2\,\pi\, {\ds{\varphi(t/\alpha)\over Q}}\Bigr)}=0,
$$
condition (\ref{eq:(4.45)}) is satisfied, and thus the asymptotic 
formula (\ref{eq:(4.59)}) holds. 
\par
Figures \ref{figure:g_2_0.1_3}$\div$\ref{figure:g_2.5_0.1_3} show the plots 
of the simulated FPT density $\tilde{g}(t)$ for 
the periodic  boundary $S(t)=S_0+B\,\sin(2\,\pi\, t/Q)$ and of the function 
$R[Z(t)]\,$ $\exp\{-\int_0^tR[Z(\tau)]\,d\tau\}$, with $R[Z(t)]$ given in 
(\ref{es:R[Z(t)]}) for various choices of parameters $S_0$, $B$ and $Q$. 
Figure \ref{figure:g_2_0.1_3} refers to  the case $S_0=2$, $B=0.1 $ and $Q=3$, 
Figure \ref{figure:g_2_0.5_3} to $S_0=2$, $B=0.5$ and $Q=3$, Figure 
\ref{figure:g_2_1_3} to  $S_0=2$, $B=1$ and $Q=3$ and  
Figure \ref{figure:g_2.5_0.1_3} to  
$S_0=2.5$, $B=0.1$ and $Q=3$.  Note that  already from rather small times,  
$\tilde{g}(t)$ is susceptible of an excellent non-homogeneous exponential 
approximation. 
%
\section{Concluding Remarks}
\setcounter{equation}{0}
Gaussian processes play an important role in numerous fields. Although many 
of their properties have been deeply analyzed, very little exists in the literature 
concerning the first-passage-time probability density function in the presence of 
constant or time-varying boundaries, which would be suitable to make predictions on 
a variety of systems evolving in the presence of some critical regions of their 
state-space.
\par
In the present paper, starting from some existing contributions to this problem area, 
for a class  of stationary Gaussian processes it has been  proved  that  a  
non-homogeneous  exponential approximation holds for the first passage time 
probability density function  in the presence of boundaries that 
either possess a horizontal asymptote or are asymptotically periodic. Furthermore, 
for a stationary Gaussian process with zero mean and with damped 
oscillatory covariance originating at $x_0=0$, extensive simulations have indicated that 
the FPT pdf  $\tilde{g}(t)$ is susceptible of an excellent non-homogeneous 
exponential approximation for either constant or periodic boundaries, even though  
these are not very distant from the initial value of the process. 
\par
We trust that such results may prove useful for the description of the time 
evolution of systems characterized, for instance, by relaxation times much 
smaller than the mean observation times, thus making particularly appropriate 
and effective the asymptotic approximation obtained in the foregoing.
\appendix
\section{Appendix I}
\setcounter{equation}{0}
Here we prove relations (\ref{eq:thmconstant10}), (\ref{eq:thmconstant7}), 
(\ref{eq:thmconstant9}) and (\ref{eq:thmconstant11ter})  of Theorem 
\ref{thmconst}.
\par
Recalling the first of (\ref{eq:componentW_n}) and  making use 
of (\ref{eq:asymconstboundary}) we have:
$$
\lim_{S_0 \to +\infty} \dot{\psi}\left( \vartheta_i\right) 
=\lim_{S_0 \to +\infty} \left[\dot{\varrho}\left( \vartheta_i\right) 
-x_0\,\dot{\gamma}\left( \vartheta_i\right)\right]=0\qquad(i=1,2,\ldots,n),
$$
where the zero limit follows from (\ref{eq:asymconstantrho}), (\ref{eq:asymR(S_0)}),   
(\ref{eq:asymconstantgamma2}) and (\ref{eq:thmconstant3bis}). This proves 
(\ref{eq:thmconstant10}).
\par
By virtue of (\ref{eq:asymR(S_0)}) and 
(\ref{eq:asymconstantgamma2}), and  recalling  (\ref{eq:thmconstant3bis}), 
from (\ref{eq:lambda})  one has 
\begin{eqnarray}
\lim_{S_0 \to +\infty} \lambda_{i+1,j+1} \left(\vartheta_1, 
\ldots, \vartheta_n \right) 
= \left\{ \begin{array}{cl}
1, & i=j=0, \ldots, n \\
- \ddot{\gamma}(0), & i=j=n+1, \ldots, 2n \\
0, & i \ne j. 
\end{array} \right.
\label{eq:thmconstant6}
\end{eqnarray}
Hence, as $S_0 \to +\infty$ the matrix $\Lambda_{2n+1} 
\left( \vartheta_1,\ldots,\vartheta_n \right)$ becomes  
diagonal with the first $n+1$ elements equal to unity and the 
last $n$ elements equal to $- \ddot{\gamma}(0)$. Therefore, 
(\ref{eq:thmconstant7}) holds.
\par 
To prove (\ref{eq:thmconstant9}), we first notice that
\begin{eqnarray}
\lim_{S_0 \to +\infty} l_{i+1,j+1} \left(\vartheta_1, 
\ldots,\vartheta_n \right)
= \left\{ \begin{array}{ll}
[-\ddot{\gamma}(0)]^n, & i=j=0, \ldots, n \\
\left[ - \ddot{\gamma}(0) \right]^{n-1}, & i=j=n+1, \ldots, 2n \\
0, & i \ne j. 
\end{array} \right. 
\label{eq:thmconstant8}
\end{eqnarray}
Therefore, by making use of (\ref{eq:thmconstant7}) and (\ref{eq:thmconstant8}), 
from the second of (\ref{eq:componentW_n}) one obtains
\begin{eqnarray*}
&&\lim_{S_0 \to +\infty} 
H\left( \vartheta_1,\ldots,\vartheta_n;\xi_1,\ldots,\xi_n \right)\\
&& \hspace*{1.0cm}= \lim_{S_0 \to +\infty} \exp\biggl\{ -{1\over 2 \left| 
\Lambda_{2n+1}(\vartheta_1,\ldots ,\vartheta_n) \right|} \sum_{i,j=1}^n 
l_{i+n+1,j+n+1}(\vartheta_1,\ldots, \vartheta_n)\,\xi_i\,\xi_j\biggr\}\\
&& \hspace*{1.0cm}=\exp \left\{ - {1\over 2\,[-\ddot{\gamma}(0)]} \sum_{i=1}^n
\xi_i^2 \right\},
\end{eqnarray*}
so that (\ref{eq:thmconstant9}) follows. 
\par
Finally, we prove (\ref{eq:thmconstant11ter}). Recalling the last of 
(\ref{eq:componentW_n}) and   (\ref{eq:R(S_0)}),  one obtains
\begin{eqnarray}
& &\hspace*{-0.6cm}{1\over \bigl[R(S_0)\bigr]^n}\,
K\left( \vartheta_1,\ldots,\vartheta_n;\xi_1,\ldots,\xi_n \right) 
={(2\pi)^n\over[- \ddot{\gamma}(0)]^{n/2}}
\label{eq:thmconstant11}\\
&&\hspace*{-0.2cm}\times  
\exp \Biggl\{ {n\,S_0^2\over 2}  -{S_0^2\over 2 \left| 
\Lambda_{2n+1}(\vartheta_1,\ldots ,\vartheta_n) \right|} 
\biggl[\,\sum_{i,j=1}^n l_{i+1,j+1}(\vartheta_1,
\ldots, \vartheta_n)\, {\psi(\vartheta_i)\over S_0}\, {\psi(\vartheta_j)\over S_0}
\nonumber\\
& &\hspace*{0.1cm}+\sum_{i,j=1}^n l_{i+1,j+n+1}(\vartheta_1,\ldots, \vartheta_n)\, 
{\psi(\vartheta_i)\over S_0^2}\,\xi_j+\sum_{i,j=1}^n l_{i+n+1,j+1}(\vartheta_1,\ldots, 
\vartheta_n)\,{\psi(\vartheta_j)\over S_0^2}\,\xi_i\biggr]\Biggr\}.\nonumber
\end{eqnarray}
We note that from the first of (\ref{eq:componentW_n}) and from 
(\ref{eq:asymconstboundary}), for $i=1,2,\ldots,n$ one has
\begin{eqnarray}
\lim_{S_0 \to +\infty}{\psi(\vartheta_i)\over S_0}=
\lim_{S_0 \to +\infty}{S(\vartheta_i)-x_0\gamma(\vartheta_i)\over S_0}
=\lim_{S_0 \to +\infty} \Bigl[{S_0+\varrho(\vartheta_i)\over S_0}
-x_0\,{\gamma(\vartheta_i)\over S_0}\Bigr]=1,
\label{eq:thmconstant11bis}
\end{eqnarray}
where the unit value follows from (\ref{eq:asymconstantrho}) , 
(\ref{eq:asymR(S_0)}), (\ref{eq:asymconstantgamma2}),  
(\ref{eq:hypthmconst}) and (\ref{eq:thmconstant3bis}). 
Taking the limit as $S_0\to +\infty$ in (\ref{eq:thmconstant11}), and 
making use of (\ref{eq:thmconstant7}), (\ref{eq:thmconstant8}) and 
(\ref{eq:thmconstant11bis}), we finally  obtain 
(\ref{eq:thmconstant11ter}).
%
\section{Appendix II}
\setcounter{equation}{0}
Here we prove relations (\ref{eq:appendix1}), (\ref{eq:appendix2}) and 
(\ref{eq:(4.54)quat}) of Theorem \ref{thmperiodic}. 
\par
By virtue of (\ref{eq:asymconstantgamma2}),  (\ref{eq:(4.50ter}) and {\it (iv)} of 
Proposition \ref{lemmaphi}, relation (\ref{eq:thmconstant6})  again follows 
from (\ref{eq:lambda}). 
Hence, as $S_0 \to +\infty$ the matrix $\Lambda_{2n+1} 
\left( \vartheta_1,\ldots,\vartheta_n \right)$ becomes  
diagonal with the first $n+1$ elements equal to unity and the 
last $n$ elements equal to $- \ddot{\gamma}(0)$, so that one immediately obtains 
(\ref{eq:appendix1}). 
\par 
The proof of (\ref{eq:appendix2})  is analogous to the proof of 
(\ref{eq:thmconstant9}) of Theorem \ref{thmconst}, taking in account 
(\ref{eq:(4.50ter}) and {\it (iv)} of Proposition \ref{lemmaphi}. 
\par
Finally, we prove (\ref{eq:(4.54)quat}). Recalling the last of 
(\ref{eq:componentW_n}),  one has
\begin{eqnarray}
& &\exp \biggl\{ {1\over2}\sum\limits_{i=1}^n
\bigl[ S_0+Z\bigl( \vartheta_i \bigr) \bigr]^2\biggr\}\,
K\left( \vartheta_1,\ldots,\vartheta_n;\xi_1,\ldots,\xi_n \right)
\label{eq:(4.54)}\\
&&= \exp \Biggl\{ {S_0^2\over2}\sum\limits_{i=1}^n
\biggl[ {S_0+Z\bigl( \vartheta_i \bigr)\over S_0} \biggr]^2  -{S_0^2\over 2 \left| 
\Lambda_{2n+1}(\vartheta_1,\ldots ,\vartheta_n) \right|}\nonumber\\ 
&&\times\biggl[\,\sum_{i,j=1}^n l_{i+1,j+1}(\vartheta_1,
\ldots, \vartheta_n)\, {\psi(\vartheta_i)\over S_0+Z\bigl( \vartheta_i \bigr)}\;
{S_0+Z\bigl( \vartheta_i \bigr)\over S_0}\; 
{\psi(\vartheta_j)\over S_0+Z\bigl( \vartheta_j \bigr)}\;
{S_0+Z\bigl( \vartheta_j \bigr)\over S_0}
\nonumber\\
& &\hspace*{0.5cm}+\sum_{i,j=1}^n l_{i+1,j+n+1}(\vartheta_1,\ldots, \vartheta_n)\, 
{\psi(\vartheta_i)\over S_0+Z\bigl( \vartheta_i\bigr)}\;
{S_0+Z\bigl( \vartheta_i\bigr)\over S_0^2}\;\;\xi_j\nonumber\\
& &\hspace*{0.5cm}+\sum_{i,j=1}^n l_{i+n+1,j+1}(\vartheta_1,\ldots, 
\vartheta_n)\,{\psi(\vartheta_j)\over S_0+Z\bigl( \vartheta_j\bigr)}\;
{S_0+Z\bigl( \vartheta_j\bigr)\over S_0^2}\;\;\xi_i\biggr]\Biggr\}.\nonumber
\end{eqnarray}
Moreover, 
\begin{eqnarray}
\lim_{S_0\to +\infty }{S_0+Z\bigl( \vartheta_i \bigr)\over S_0}=1
\qquad(i=1,2,\ldots,n),
\label{eq:(4.54)bis}
\end{eqnarray}
since $Z(t)$ is a bounded function independent of $S_0$. Furthermore, we note that 
from the first of (\ref{eq:componentW_n}) and from 
(\ref{eq:asymperiodicboundary}), for $i=1,2,\ldots,n$ one obtains
\begin{eqnarray}
&&\lim_{S_0\to +\infty }{\psi(\vartheta_i)\over S_0+Z\bigl( \vartheta_i\bigr)}
=\lim_{S_0\to +\infty }{S(\vartheta_i)-x_0\,\gamma(\vartheta_i)
\over S_0+Z\bigl( \vartheta_i\bigr)}\nonumber\\
&& \hspace*{3cm}=\lim_{S_0\to +\infty }\Bigl[{S_0+\varrho(\vartheta_i)
\over S_0+Z\bigl( \vartheta_i\bigr)}
-x_0\;{\gamma(\vartheta_i)\over S_0+Z\bigl( \vartheta_i\bigr)}\Bigr]=1,
\label{eq:(4.54)ter}
\end{eqnarray}
where the last equality follows by using condition {\it (iii)} of 
Proposition \ref{lemmaphi} and by 
recalling  (\ref{eq:asymconstantgamma2}), (\ref{eq:Z(t)}) and  (\ref{eq:(4.45)}). 
Taking the limit as $S_0\to +\infty$ in (\ref{eq:(4.54)}), and 
making use of (\ref{eq:thmconstant7}), (\ref{eq:thmconstant8}), 
(\ref{eq:(4.54)bis}) and (\ref{eq:(4.54)ter}),  one  thus obtains 
(\ref{eq:(4.54)quat}).
%
%

\end{document}